\documentclass[12pt]{article}
\usepackage{tikz}
\usetikzlibrary{arrows}
\usepackage[framemethod=tikz]{mdframed}
\usepackage{wrapfig}
\usepackage{amsmath,amssymb}
\usepackage{type1cm}
\usepackage{amsthm}
\usepackage{mathrsfs}
\usepackage{mathtools}
\usepackage{enumerate}
\usepackage[all]{xy} 
\usepackage[vcentermath,enableskew]{youngtab}
\usepackage{ytableau}
\usepackage{diagbox}
\AtBeginDocument{%
   \def\MR#1{}
}

\DeclareMathOperator{\GL}{GL}

\DeclareMathOperator{\ch}{ch}

\DeclareMathOperator{\Adm}{Adm}

\DeclareMathOperator{\Irr}{Irr}
\DeclareMathOperator{\pr}{pr}
\DeclareMathOperator{\wt}{wt}
\DeclareMathOperator{\red}{red}
\DeclareMathOperator{\pfn}{pfn}
\DeclareMathOperator{\FE}{FE}
\DeclareMathOperator{\op}{op}
\DeclareMathOperator{\dom}{dom}
\DeclareMathOperator{\en}{end}

\DeclareMathOperator{\supp}{supp}
\DeclareMathOperator{\de}{def}
\DeclareMathOperator{\LP}{LP}
\DeclareMathOperator{\cyc}{cyc}

\newcommand\F{\mathbb{F}}

\newcommand\aFq{\overline{\mathbb F}_q}

\newcommand\cA{\mathcal A}
\newcommand\cB{\mathcal B}
\newcommand\cO{\mathcal O}
\newcommand\cI{\mathcal I}
\newcommand\cV{\mathcal V}
\newcommand\cG{\mathcal G}

\newcommand\cT{\mathcal T}

\newcommand\Gm{\mathbb G_m}

\newcommand\N{\mathbb N}
\newcommand\Q{\mathbb Q}

\newcommand\A{\mathbb A}
\newcommand\B{\mathbb B}
\newcommand\G{\mathbb G}
\newcommand\J{\mathbb J}
\newcommand\q{\bold q}
\newcommand\hG{\widehat G}

\newcommand\Z{\mathbb Z}
\newcommand\tW{\tilde W}

\newcommand\te{\tilde e}
\newcommand\tf{\tilde f}
\newcommand\SW{{^S\tilde W}}
\newcommand\SAdm{{^S\mathrm{Adm}}}

\newcommand\tS{\tilde S}
\newcommand\inv{\mathrm{inv}}
\newcommand\bl{\bullet}
\newcommand\tp{\mathrm{top}}
\newcommand\ld{\lambda}
\newcommand\up{\upsilon}
\newcommand\unp{\underline p}
\newcommand\ve{\varepsilon}
\newcommand\vph{\varphi}
\newcommand\vp{\varpi}
\newcommand\Y{X_*(T)}

\newcommand\bb{\bold b}
\newcommand\mn{\lfloor \frac{m}{n} \rfloor}
\newcommand\la{\langle}
\newcommand\ra{\rangle}
\newcommand\pc{\preceq}

\theoremstyle{definition}
\newtheorem{theo}{Theorem}[section]
\newtheorem{prop}[theo]{Proposition}
\newtheorem{defi}[theo]{Definition}
\newtheorem{lemm}[theo]{Lemma}
\newtheorem{coro}[theo]{Corollary}
\newtheorem{exam}[theo]{Example}
\newtheorem{rema}[theo]{Remark}

\newtheorem{thm}{Theorem}[section]

\begin{document}
\title{The Ekedahl-Oort Stratification and the Semi-Module Stratification}
\author{Ryosuke Shimada}
\date{}
\maketitle

\begin{abstract}
In this paper we compare the $\J$-stratification (or the semi-module stratification) and the Ekedahl-Oort stratification of affine Deligne-Lusztig varieties in the superbasic case.
In particular, we classify the cases where the $\J$-stratification gives a refinement of the Ekedahl-Oort stratification, which include many interesting cases such that the affine Deligne-Lusztig variety admits a simple geometric structure.
\end{abstract}

\section{Introduction}
\label{introduction}
The affine Deligne-Lusztig variety was introduced by Rapoport in \cite{Rapoport}, which plays an important role in understanding geometric and arithmetic properties of Shimura varieties.
The uniformization theorem by Rapoport and Zink \cite{RZ} allows us to describe the Newton strata of Shimura varieties in terms of Rapoport-Zink spaces, whose underlying spaces are special cases of affine Deligne-Lusztig varieties.

Let $F$ be a non-archimedean local field with finite field $\F_q$ of prime characteristic $p$, and let $L$ be the completion of the maximal unramified extension of $F$.
Let $\sigma$ denote the Frobenius automorphism of $L/F$.
Further, we write $\cO$ (resp.\ $\cO_F$) for the valuation ring $L$ (resp.\ $F$).
Finally, we denote by $\vp$ a uniformizer of $F$ (and $L$) and by $v_L$ the valuation of $L$ such that $v_L(\vp)=1$.

Let $G$ be an unramified connected reductive group over $\cO_F$.
Let $B\subset G$ be a Borel subgroup and $T\subset B$ a maximal torus in $B$, both defined over $\cO_F$.
For $\mu,\mu'\in X_*(T)$ (resp.\ $X_*(T)_{\Q}$), we write $\mu'\pc \mu$ if $\mu-\mu'$ is a non-negative integral (resp.\ rational) linear combination of positive coroots.
For a cocharacter $\mu\in X_*(T)$, let $\vp^{\mu}$ be the image of $\vp\in \mathbb G_m(F)$ under the homomorphism $\mu\colon\mathbb G_m\rightarrow T$.

Set $K=G(\cO)$.
We fix a dominant cocharacter $\mu\in X_*(T)_+$ and $b\in G(L)$.
Then the affine Deligne-Lusztig variety $X_{\mu}(b)$ is the locally closed reduced $\aFq$-subscheme of the affine Grassmannian $\cG r=G(L)/K$ defined as
$$X_{\mu}(b)=\{xK\in \cG r\mid x^{-1}b\sigma(x)\in K\vp^{\mu}K\}.$$
The closed affine Deligne-Lusztig variety is the closed reduced $\aFq$-subscheme of $\cG r$ defined as
$$X_{\pc\mu}(b)=\bigcup_{\mu'\pc \mu}X_{\mu'}(b).$$
Both $X_{\mu}(b)$ and $X_{\pc\mu}(b)$ are locally of finite type in the equal characteristic case and locally perfectly of finite type in the mixed characteristic case (cf.\ \cite[Corollary 6.5]{HV}, \cite[Lemma 1.1]{HV2}).
Finally, the affine Deligne-Lusztig varieties $X_{\mu}(b)$ and $X_{\pc\mu}(b)$ carry a natural action (by left multiplication) by the $\sigma$-centralizer of $b$
$$J_b(F)=\{g\in G(L)\mid g^{-1}b\sigma(g)=b\}.$$

The geometric properties of affine Deligne-Lusztig varieties have been studied by many people.
For example, the non-emptiness criterion and the dimension formula are already known for the affine Deligne-Lusztig varieties in the affine Grassmannian (see \cite{Gashi}, \cite{Viehmann} and \cite{Hamacher}).
Let $B(G)$ denote the set of $\sigma$-conjugacy classes of $G(L)$. 
Thanks to Kottwitz \cite{Kottwitz2}, a $\sigma$-conjugacy class $[b]\in B(G)$ is uniquely determined by two invariants: the Kottwitz point $\kappa(b)\in \pi_1(G)/((1-\sigma)\pi_1(G))$ and the Newton point $\nu_b\in X_*(T)_{\Q,+}$.
Set $B(G,\mu)=\{[b]\in B(G)\mid \kappa(b)=\kappa(\vp^\mu), \nu_b\pc \mu^{\diamond}\}$, where $\mu^{\diamond}$ denotes the $\sigma$-average of $\mu$.
Then $X_\mu(b)\neq \emptyset$ if and only if $[b]\in B(G,\mu)$.
If this is the case, then we have
$$\dim X_\mu(b)=\la\rho, \mu-\nu_b\ra-\frac{1}{2}\de(b),$$
where $\rho$ is the half sum of positive roots and $\de(b)$ is the defect of $b$.
Moreover, the parametrization problem of the set of irreducible components $\Irr X_\mu(b)$ is also known.
Let $\hG$ be the Langlands dual of $G$ defined over $\overline{\Q}_l$ with $l\neq p$.
Surprisingly, there exists a natural bijection between $J_b(F)\backslash\Irr X_\mu(b)$ and a certain weight space of the crystal basis $\B_\mu$ of the irreducible $\hG$-module $V_\mu$ of highest weight $\mu$.
This is conjectured by Chen and Zhu, and proved in general by Nie \cite{Nie}.

Besides them, it is known that in certain cases, the closed affine Deligne-Lusztig variety admits a simple description.
In \cite{GH}, \cite{GHN} and \cite{GHN3}, a notion of ``Coxeter type'' was introduced by G\"{o}rtz, He and Nie.
They proved that if $(G, \mu)$ is of Coxeter type and if $b$ is the unique basic element in $B(G,\mu)$, then $X_{\pc \mu}(b)$ is naturally a union of classical Deligne-Lusztig varieties (in fact, they studied the cases with arbitrary parahoric level).
This stratification is the so-called Bruhat-Tits stratification, a stratification indexed in terms of the Bruhat-Tits building of $J_b(F)$, see \cite[\S 2.4]{GHN3}.
These simple descriptions of closed affine Deligne-Lusztig varieties have been applied to number theory especially when $(G,\mu)$ corresponds to a Shimura datum (cf.\ \cite[\S1]{GHN3}).
For example, the cases of Coxeter type include the case for certain unitary groups of signature $(1, n-1)$ studied in \cite{VW} by Vollaard and Wedhorn, which has been used in the Kudla-Rapoport program \cite{KR11}.

To give a conceptual way to explain the relationship between the geometry of affine Deligne-Lusztig varieties and the Bruhat-Tits building of $J_b(F)$ indicated by above examples, Chen and Viehmann \cite{CV} introduced the $\J$-stratification, where $\J$ stands for $J_b(F)$.
The $\J$-strata are locally closed subsets of $\cG r$.
By intersecting each $\J$-stratum with $X_{\pc\mu}(b)$, we obtain the $\J$-stratification of $X_{\pc\mu}(b)$ (see \S\ref{J-str} for details).
In \cite{Gortz2}, G\"{o}rtz showed that the Bruhat-Tits stratification coincides with the $\J$-stratification.
In fact the Bruhat-Tits stratification is a refinement of the Ekedahl-Oort stratification (see \S\ref{ADLV} for the latter).
So the $\J$-stratification is also a refinement of the Ekedahl-Oort stratification when $(G,\mu)$ is of Coxeter type.
This does not hold in general even if $\mu$ is minuscule.
See \cite[Example 4.1]{CV} for a counterexample in the case $G=\GL_9$.
Therefore the cases when $\J$-stratification is a refinement of the Ekedahl-Oort stratification should be special cases, which are of particular interest.

Usually it seems very difficult to study the $\J$-stratification.
However, in the case that $G=\GL_n$ and $b$ is superbasic (i.e., $\kappa(b)\in \Z$ is coprime to $n$), the $\J$-stratification coincides with a stratification by semi-modules (\cite[Proposition 3.4]{CV}).
The notion of semi-modules was first considered by de Jong and Oort \cite{JO} (see \S\ref{exsemi}) for minuscule cocharacters.
Later Viehmann \cite{Viehmann} introduced a notion of extended semi-modules for arbitrary cocharacters, which generalizes the notion of semi-modules.
It played a crucial role to prove the dimension formula (for split groups) and the Chen-Zhu conjecture mentioned above.
This is because for these problems, we can reduce the general case to the case that $G=\GL_n$ and $b$ is superbasic.

The aim of this paper is to compare the Ekedahl-Oort stratification and the semi-module stratification (for $G=\GL_n$).
To state the main results, we need some notation.
Let $W_0$ be the (finite) Weyl group of $T$ in $G$ and let $\tW$ be the Iwahori-Weyl group of $T$ in $G$.
Then $\tW=X_*(T)\rtimes W_0$.
We denote the projection $\tW\rightarrow W_0$ by $p$.
For $\mu\in \Y_+$, we denote by $\Adm(\mu)$ the admissible subset of $\tW$.
Let $\SAdm(\mu)$ be a certain subset of $\Adm(\mu)$, which is the index set of the Ekedahl-Oort stratification of $X_{\pc\mu}(\tau_\mu)$ (see \S\ref{ADLV}).
We fix (a representative in $G(L)$ of a) length $0$ element $\tau_\mu\in \tW$ whose $\sigma$-conjugacy class in $G(L)$ is the unique basic element in $B(G,\mu)$.
Finally, let $\LP(w)\subseteq W_0$ be the length positive elements for $w$ (see \S\ref{LP}).
\begin{thm}[See Theorem \ref{main theo}]
\label{main thm}
Let $G=\GL_n$ and let $\mu\in \Y_+$.
Assume that $\tau_\mu$ is superbasic.
Then the following assertions are equivalent.
\begin{enumerate}[(i)]
\item The $\J$-stratification (or the semi-module stratification) of $X_{\pc\mu}(\tau_\mu)(\neq \emptyset)$ gives a refinement of the Ekedahl-Oort stratification.
\item For any $w\in {^S{\Adm}}(\mu)$ whose corresponding Ekedahl-Oort stratum is non-empty, there exists $v\in \LP(w)$ such that $v^{-1}p(w)v$ is a Coxeter element.
\item The cocharacter $\mu$ has one of the following forms modulo $\Z\omega_n$:
\begin{align*}
&\omega_1,\quad \omega_{n-1},\ &(n\geq 1),\\
&\omega_2,\quad 2\omega_1,\quad \omega_{n-2},\quad 2\omega_{n-1},\ &(\text{odd}\ n\geq 3),\\
&\omega_2+\omega_{n-1},\quad 2\omega_1+\omega_{n-1}\quad \omega_1+\omega_{n-2},\quad\omega_1+2\omega_{n-1},\ &(n\geq 3),\\
&\omega_3,\quad\omega_{n-3},\ &(n=7,8),\\
&3\omega_1,\quad 3\omega_{n-1},\ &(n=4,5),\\
&\omega_1+\omega_2,\quad\omega_3+\omega_4,\  &(n=5),\\
&4\omega_1,\quad \omega_1+3\omega_2,\quad 4\omega_2,\quad 3\omega_1+\omega_2, &(n=3),\\
&m\omega_1\ \text{with $m$ odd,} &(n=2).
\end{align*}
\end{enumerate}
Here $\omega_k$ denotes the cocharacter of the form $(1,\ldots,1,0,\ldots,0)$ in which $1$ is repeated $k$ times.
\end{thm}

See \S\ref{J-str} for the reason why we choose $\tau_\mu$.
In fact, this choice is the reasonable one suggested in \cite[Remark 2.1]{CV}, which is unique in this case.
Also we can deduce the geometric structure of each $\J$-stratum.
In many cases, it is universally homeomorphic to (the perfection of) an affine space (see Theorem \ref{main theo} and Remark \ref{main rema}).

Although the cocharacters $\omega_1$ and $\omega_{n-1}$ are of Coxeter type for any $n$, the cocharacters $2\omega_1$ and $\omega_2$ are of Coxeter type only when $n=2$ and $n=4$ respectively (cf.\ \cite[Theorem 1.4]{GHN3}).
In Theorem \ref{main thm}, these two cocharacters are no longer exceptional cases.
Note also that the condition (ii) works for general $G$.
It would be interesting to study this condition in general.

Cyclic semi-modules are certain simple elements in the set of extended semi-modules.
It is easy to see that if there exists a non-cyclic semi-module for $\mu$, then the semi-module stratification of $X_\mu(\tau_\mu)$ never gives a refinement of the Ekedahl-Oort stratification (Corollary \ref{never refinement}).
Along the way of proving Theorem \ref{main thm}, we also prove the following classification theorem, which ensures that there exists a non-cyclic semi-module in many cases.

\begin{thm}[See Theorem \ref{cyclic classification}]
\label{main thm2}
Every top extended semi-module (the semi-module whose corresponding stratum is top-dimensional) for $\mu$ is cyclic if and only if $\mu$ has one of the following forms modulo $\Z\omega_n$:
\begin{enumerate}[(i)]
\item $\omega_i$ with $1\le i\le n-1$ such that $i$ is coprime to $n$.
\item $\omega_1+\omega_i$ or $\omega_{n-1}+\omega_{n-i}$  with $1\le i\le n-1$ such that $i+1$ is coprime to $n$.
\item $(nr+i)\omega_1$ or $(nr+i)\omega_{n-1}$ with $r\geq 0$ and $1\le i\le n-1$ such that $i$ is coprime to $n$.
\item $(nr+i-j)\omega_1+\omega_j$ or $(nr+i-j)\omega_{n-1}+\omega_{n-j}$ with $r\geq 1$, $2\le j\le n-1$ and $1\le i\le n-1$ such that $i$ is coprime to $n$.
\end{enumerate}
\end{thm}

The paper is organized as follows.
In \S\ref{preliminaries} we introduce the affine Deligne-Lusztig variety and stratifications of it.
We also recall the length positive elements and the non-emptiness criterion of the affine Deligne-Lusztig variety in the affine flag variety.
In \S\ref{semi-modules} and \S\ref{crystal bases}, we recollect known results on semi-modules and crystal bases respectively.
Also in \S\ref{crystal bases}, we prove Theorem \ref{main thm2} using combinatorics on Young tableaux.
The key ingredient here is the explicit construction of top extended semi-modules from crystal bases via the natural map in the Chen-Zhu conjecture.
This is established in \cite{Shimada3} by the author.
In \S\ref{the semi-module stratification} and \S\ref{the Ekedahl-Oort stratification}, we examine the semi-module stratification and the Ekedahl-Oort stratification respectively by an explicit calculation of semi-modules and elements in $\SAdm(\mu)$.
In particular, using the non-emptiness criterion mentioned above, we show that Theorem \ref{main thm} (ii) does not hold for many $\mu$.
Finally in \S7 we prove the main theorem, combining Theorem \ref{main thm2} and the results in \S\ref{the semi-module stratification} and \S\ref{the Ekedahl-Oort stratification}.

\textbf{Acknowledgments:}
The author is grateful to his advisor Yoichi Mieda for his constant support and encouragement.
This work was supported by the WINGS-FMSP program at the Graduate School of Mathematical Science, the University of Tokyo. 
This work was also supported by JSPS KAKENHI Grant number JP21J22427.

\section{Preliminaries}
\label{preliminaries}
Keep the notations in \S\ref{introduction}.
\subsection{Notation}
\label{notation}
Let $\Phi=\Phi(G,T)$ denote the set of roots of $T$ in $G$.
We denote by $\Phi_+$ (resp.\ $\Phi_-$) the set of positive (resp.\ negative) roots distinguished by $B$.
Let $\Delta$ be the set of simple roots and $\Delta^\vee$ be the corresponding set of simple coroots.
Let $X_*(T)$ be the set of cocharacters, and let $X_*(T)_+$ be the set of dominant cocharacters.

The Iwahori-Weyl group $\tW$ is defined as the quotient $N_{G(L)}T(L)/T(\cO)$.
This can be identified with the semi-direct product $W_0\ltimes X_{*}(T)$, where $W_0$ is the finite Weyl group of $G$.
We denote the projection $\tW\rightarrow W_0$ by $p$.
Let $S\subset W_0$ denote the subset of simple reflections, and let $\tS\subset \tW$ denote the subset of simple affine reflections.
We often identify $\Delta$ and $S$.
The affine Weyl group $W_a$ is the subgroup of $\tW$ generated by $\tS$.
Then we can write the Iwahori-Weyl group as a semi-direct product $\tW=W_a\rtimes \Omega$, where $\Omega\subset \tW$ is the subgroup of length $0$ elements.
Moreover, $(W_a, \tS)$ is a Coxeter system.
We denote by $\le$ the Bruhat order on $\tW$.
For any $J\subseteq \tS$, let $^J\tW$ be the set of minimal length elements for the cosets in $W_J\backslash \tW$, where $W_J$ denotes the subgroup of $\tW$ generated by $J$.
We also have a length function $\ell\colon \tW\rightarrow \Z_{\geq 0}$ given as
$$\ell(w_0\vp^{\lambda})=\sum_{\alpha\in \Phi_+, w_0\alpha\in \Phi_-}|\langle \alpha, \lambda\rangle+1|+\sum_{\alpha\in \Phi_+, w_0\alpha\in \Phi_+}|\langle \alpha, \lambda\rangle|,$$
where $w_0\in W_0$ and $\lambda\in \Y$.

For $w\in W_a$, we denote by $\supp(w)\subseteq \tS$ the set of simple affine reflections occurring in every (equivalently, some) reduced expression of $w$.
Note that $\tau\in \Omega$ acts on $\tS$ by conjugation.
We define the $\sigma$-support $\supp_\sigma(w\tau)$ of $w\tau$ as the smallest $\tau\sigma$-stable subset of $\tS$.
We call an element $w\tau\in W_a\tau$ a $\sigma$-Coxeter element if exactly one simple reflection from each $\tau\sigma$-orbit on $\supp_\sigma(w\tau)$ occurs in every (equivalently, any) reduced expression of $w$.

For $w,w'\in \tW$ and $s\in \tS$, we write $w\xrightarrow{s}_\sigma w'$ if $w'=sw\sigma(s)$ and $\ell(w')\le \ell(w)$.
We write $w\rightarrow_\sigma w'$ if there is a sequence $w=w_0,w_1,\ldots, w_k=w'$ of elements in $\tW$ such that for any $i$, $w_{i-1}\xrightarrow{s_i}_\sigma w_i$ for some $s_i\in S$.
If $w\rightarrow_\sigma w'$ and $w'\rightarrow_\sigma w$, we write $w\approx_\sigma w'$.

For $\alpha\in \Phi$, let $U_\alpha\subseteq G$ denote the corresponding root subgroup.
We set $$I=T(\cO)\prod_{\alpha\in \Phi_+}U_{\alpha}(\vp\cO)\prod_{\beta\in \Phi_-}U_{\beta}(\cO)\subseteq G(L),$$
which is called the standard Iwahori subgroup associated to the triple $T\subset B\subset G$.

In the case $G=\GL_n$, we will use the following description.
Let $\chi_{ij}$ be the character $T\rightarrow \Gm$ defined by $\mathrm{diag}(t_1,t_2,\ldots, t_n)\mapsto t_i{t_j}^{-1}$.
Then we have $\Phi=\{\chi_{ij}\mid i\neq j\}$, $\Phi_+=\{\chi_{ij}\mid i< j\}$, $\Phi_-=\{\chi_{ij}\mid i> j\}$ and $\Delta=\{\chi_{i,i+1}\mid 1\le i <n\}$.
Through the isomorphism $X_*(T)\cong \Z^n$, ${X_*(T)}_+$ can be identified with the set $\{(m_1,\cdots, m_n)\in \Z^n|m_1\geq \cdots \geq m_n\}$.
Let us write $s_1=(1\ 2), s_2=(2\ 3), \ldots, s_{n-1}=(n-1\ n)$.
Set $s_0=\vp^{\chi_{1,n}^{\vee}}(1\ n)$, where $\chi_{1,n}$ is the unique highest root.
Then $S=\{s_1,s_2,\ldots, s_{n-1}\}$ and $\tS=S\cup\{s_0\}$.
The Iwahori subgroup $I\subset K$ is the inverse image of the lower triangular matrices under the projection $G(\cO)\rightarrow G(\aFq),\  \vp\mapsto 0$.
Set $\tau={\begin{pmatrix}
0 & \vp \\
1_{n-1} & 0\\
\end{pmatrix}}$.
We often regard $\tau$ as an element of $\tW$, which is a generator of $\Omega\cong \Z$.
Note that $b\in \GL_n(L)$ is superbasic if and only if $[b]=[\tau^m]$ in $B(\GL_n)$ for some $m$ coprime to $n$.

\subsection{Affine Deligne-Lusztig Varieties}
\label{ADLV}
For $w\in \tW$ and $b\in G(L)$, the affine Deligne-Lusztig variety $X_w(b)$ in the affine flag variety $G(L)/I$ is defined as
$$X_w(b)=\{xI\in G(L)/I\mid x^{-1}b\sigma(x)\in IwI\}.$$
For $\mu\in \Y_+$ and $b\in G(L)$, the affine Deligne-Lusztig variety $X_{\mu}(b)$ in the affine Grassmannian $\cG r=G(L)/K$ is defined as
$$X_{\mu}(b)=\{xK\in \cG r\mid x^{-1}b\sigma(x)\in K\vp^{\mu}K\}.$$
The closed affine Deligne-Lusztig variety is the closed reduced $\aFq$-subscheme of $\cG r$ defined as
$$X_{\pc\mu}(b)=\bigcup_{\mu'\pc \mu}X_{\mu'}(b).$$
Left multiplication by $g^{-1}\in G(L)$ induces an isomorphism between $X_\mu(b)$ and $X_\mu(g^{-1}b\sigma(g))$.
Thus the isomorphism class of the affine Deligne-Lusztig variety only depends on the $\sigma$-conjugacy class of $b$.
Moreover, we have $X_\mu(b)=X_{\mu+\ld}(\vp^{\ld}b)$ for each central $\ld \in \Y$.

The admissible subset of $\tW$ associated to $\mu$ is defined as
$$\Adm(\mu)=\{w\in \tW\mid w\le \vp^{w_0\mu}\ \text{for some}\ w_0\in W_0\}.$$
Set $\SAdm(\mu)=\Adm(\mu)\cap \SW$.
Then, by \cite[Theorem 3.2.1]{GH} (see also \cite[\S2.5]{GHR}), we have
$$X_{\pc\mu}(b)=\bigsqcup_{w\in\SAdm(\mu)}\pi(X_w(b)),$$
where $\pi\colon G(L)/I\rightarrow G(L)/K$ is the projection.
This is the so-called Ekedahl-Oort stratification.

For any $w\in \SW$, set 
$$Z(w)\coloneqq\{w_0\in W_0\mid w_0w=ww_0\}.$$ 
\begin{lemm}
\label{injective}
Let $\vp^\mu y\in \SW$ with $\mu$ dominant and $y\in W_0$.
Assume that $Z(\vp^\mu y)=\{1\}$.
Then the projection map $\pi\colon X_{\vp^\mu y}(b)\rightarrow X_\mu(b)$ is injective.
\end{lemm}
\begin{proof}
The proof is similar to \cite[Lemma 5.4]{HNY}.
We may assume that $X_{\vp^\mu y}(b)\neq \emptyset$.
Let $gI,g'I\in X_{\vp^\mu y}(b)$ such that $\pi(gI)=\pi(g'I)$.
Then $g'^{-1}g\in K$ and hence $g'^{-1}g\in IxI$ for some $x\in W_0$.
Since $(g'^{-1}g)(g^{-1}b\sigma(g))=(g'^{-1}b\sigma(g'))(\sigma(g'^{-1}g))$, we have $(IxI)(I\vp^\mu yI)\cap (I\vp^\mu yI)(IxI)\neq \emptyset$.
Note that $(IxI)(I\vp^\mu yI)=Ix \vp^\mu yI$ because $\vp^\mu y\in \SW$.
This implies that $x\vp^\mu y=\vp^\mu y x$.
By our assumption, we must have $x=1$ and hence $g'^{-1}g\in I$ as desired.
\end{proof}

\begin{exam}
Let $G=\GL_n$ and $b=\tau^m$ with $m$ coprime to $n$,
and let $\vp^\mu y\in \SW$ with $\mu$ dominant and $y\in W_0$.
If $\mu$ is minuscule, then it is easy to check that $\vp^\mu y=\tau^m y'$ for some $y'\in W_0$.
Since $m$ is coprime to $n$, we conclude that $Z(\vp^\mu y)=\{1\}$.
If $y$ is $n$-cycle and $\{s_1,s_{n-1}\}\nsubseteq Z(\vp^\mu)$, then we also have $Z(\vp^\mu y)=\{1\}$.
Indeed, for any $x\in W_0$, $x\vp^\mu y=\vp^\mu y x$ implies that $xyx^{-1}=y$ and $x\in Z(\vp^\mu)$.
Thus $x=y^k$ for some $0\le k\le n-1$ and $y^k\mu=\mu$.
Since $\{s_1,s_{n-1}\}\nsubseteq Z(\vp^\mu)$, we must have $k=0$.
\end{exam}

\subsection{Deligne-Lusztig Reduction Method}
\label{DL method}
The following Deligne-Lusztig reduction method was established in \cite[Corollary 2.5.3]{GH3}.
\begin{prop}
\label{DL method prop}
Let $w\in \tW$ and let $s\in \tS$ be a simple affine reflection.
If $\ch(F)>0$, then the following two statements hold for any $b\in G(L)$.
\begin{enumerate}[(i)]
\item If $\ell(sw\sigma(s))=\ell(w)$, then there exists a $J_b(F)$-equivariant universal homeomorphism $X_w(b)\rightarrow X_{sw\sigma(s)}(b)$.
\item If $\ell(sw\sigma(s))=\ell(w)-2$, then there exists a decomposition $X_w(b)=X_1\sqcup X_2$ such that
\begin{itemize}
\item $X_1$ is open and there exists a $J_b(F)$-equivariant morphism $X_1\rightarrow X_{sw}(b)$, which is  the composition of a Zariski-locally trivial $\G_m$-bundle and a universal homeomorphism. 
\item $X_2$ is closed and there exists a $J_b(F)$-equivariant morphism $X_2\rightarrow X_{sw\sigma(s)}(b)$, which is the composition of a Zariski-locally trivial $\A^1$-bundle and a universal homeomorphism. 
\end{itemize}
If $\ch(F)=0$, then the above statements still hold by replacing $\A^1$ and $\G_m$ by $\A^{1,\pfn}$ and $\G_m^{\pfn}$ respectively.
\end{enumerate}
\end{prop}

The following result is proved in \cite[Theorem 2.10]{HN2}, which allows us to reduce the study of $X_w(b)$ for any $w$, via the Deligne-Lusztig reduction method, to the study of $X_w(b)$ for $w$ of minimal length in its $\sigma$-conjugacy class.
\begin{theo}
\label{minimal}
For each $w\in \tW$, there exists an element $w'$ which is of minimal length inside its $\sigma$-conjugacy class such that $w\rightarrow_\sigma w'$.
\end{theo}

Following \cite[\S 3.4]{HNY}, we construct the reduction trees for $w$ by induction on $\ell(w)$.

The vertices of the trees are the elements of $\tW$.
We write $x\rightharpoonup y$ if $x,y\in \tW$ and there exists $x'\in \tW$ and $s\in \tS$ such that $x\approx_\sigma x'$, $\ell(sx'\sigma(s))=\ell(x')-2$ and $y\in \{sx', sx'\sigma(s)\}$.
These are (oriented) edges of the trees.

If $w$ is of minimal length in its $\sigma$-conjugacy class of $\tW$, then the reduction tree for $w$ consists of a single vertex $w$ and no edges.
Assume that $w$ is not of minimal length and that a reduction tree is given for any $z\in \tW$ with $\ell(z)<\ell(w)$.
By Theorem \ref{minimal}, there exist $w'$ and $s\in \tS$ with $w\approx_\sigma w'$ and $\ell(sw'\sigma(s))=\ell(w')-2$.
Then a reduction tree of $w$ consists of the given reduction trees of $sw'$ and $sw'\sigma(s)$ and the edges $w\rightharpoonup sw'$ and $w\rightharpoonup sw'\sigma(s)$. 

Let $\cT$ be a reduction tree of $w$.
An end point of $\cT$ is a vertex in $\cT$ of minimal length.
A reduction path in $\cT$ is a path $\unp\colon w \rightharpoonup w_1\rightharpoonup\cdots \rightharpoonup w_n$, where $w_n$ is an end point of $\cT$.
Set $\en(\unp)=w_n$.
We say that $x\rightharpoonup y$ is of type I (resp.\ II) if $\ell(x)-\ell(y)=1$ (resp.\ $\ell(x)-\ell(y)=2$).
For any reduction path $\unp$, we denote by $\ell_{I}(\unp)$ (resp.\ $\ell_{II}(\unp)$) the number of type I (resp.\ II) edges in $\unp$.
We write $X_{\unp}$ a locally closed subscheme of $X_w(b)$ which is $J_b(F)$-equivariant universally homeomorphic to an iterated fibration of type $(\ell_{I}(\unp),\ell_{II}(\unp))$ over $X_{\en(\unp)}(b)$.

Let $B(\tW,\sigma)$ be the set of $\sigma$-conjugacy classes in $\tW$.
Let $\Psi:B(\tW,\sigma)\rightarrow B(G)$ be the map sending $[w]\in B(\tW,\sigma)$ to $[\dot w]\in B(G)$.
It is known that this map is well-defined and surjective, see \cite[Theorem 3.7]{He14}.
By \cite[Proposition 3.9]{HNY}, we have the following description of $X_w(b)$.
\begin{prop}
\label{decomposition}
Let $w\in \tW$ and $\cT$ be a reduction tree of $w$.
For any $b\in G(L)$, there exists a decomposition
$$X_w(b)=\bigsqcup_{\substack{\unp\ \text{is a reduction path in}\ \cT;\\ \Psi(\en(\unp))=[b]}}X_{\unp}.$$
\end{prop}

In the case that $G=\GL_n$ and $b=\tau^m$ with $m$ coprime to $n$, we can count the number of top irreducible components and rational points of $X_w(b)^0=\{gI\in X_w(b)\mid \kappa(g)=v_L(\det(g))=0\}$ using the reduction tree for $w$.
By \cite[Proposition 3.5]{HN2}, the $\sigma$-conjugacy class of $\tau^m$ in $\tW$ is the unique element in $B(\tW,\sigma)$ which maps to $[\tau^m]\in B(G)$ under $\Psi$.
Note also that $\tau^m$ is the unique minimal length element in its $\sigma$-conjugacy class.
We define a polynomial as
$$F_{w,b}\coloneqq\sum_{\unp}(\q-1)^{\ell_I(\unp)}\q^{\ell_{II}(\unp)}\in \N[\q-1],$$
where $\unp$ runs over all the reduction paths in $\cT$ with $\en(\unp)=\tau^m$.

\begin{prop}
\label{rational point}
Assume that $G=\GL_n$ and $b=\tau^m$ with $m$ coprime to $n$.
Let $w\in \tW$ and let $\cT$ be a reduction tree of $w$.
Then the number of top irreducible components of $X_w(b)^0$ is equal to the leading coefficient of $F_{w,b}$.
Moreover, we have
$$|X_w(b)^{0,\sigma}|=F_{w,b}|_{\q=q}.$$
\end{prop}
\begin{proof}
Note that each $J_b(F)$-orbit of an irreducible component of $X_w(b)$ can be represented by an irreducible component of $X_w(b)^0$.
Moreover, it is known that the stabilizer in $J_b(F)$ is a parahoric subgroup (cf.\ \cite[Proposition 3.1.4]{ZZ}), i.e., $J_b(F)\cap I=\{g\in J_b(F)\mid \kappa(g)=0\}$.
Then the statement follows from \cite[Theorem 3.4 \& Proposition 3.5]{HNY} and \cite[Corollary 4.4]{HZ20}.
\end{proof}

\begin{rema}
The polynomials $F_{w,b}$ are called {\it class polynomials}.
However, the definition above is an ad hoc one.
See \cite[\S3]{HNY} for the definition in general and the connection to reduction trees.
\end{rema}

\subsection{The $\J$-stratification}
\label{J-str}
For any $g,h\in G(L)$, let $\inv(g,h)$ denote the relative position, i.e., the unique dominant cocharacter such that $g^{-1}h\in K\vp^{\inv(g,h)} K$.
By definition, two elements $gK,hK\in G(L)/K$ lie in the same $\J$-stratum if and only if for all $j\in J_b(F)$, $\inv(j,g)=\inv(j,h)$.
Clearly, this does not depend on the choice of $g,h$.
By \cite[Proposition 2.11]{CV}, the $\J$-strata are locally closed in $\cG r$.
By intersecting each $\J$-stratum with $X_\mu(b)$ (resp.\ $X_{\pc\mu}(b)$), we obtain the $\J$-stratification of $X_\mu(b)$ (resp.\ $X_{\pc\mu}(b)$).

As explained in \cite[Remark 2.1]{CV}, the $\J$-stratification heavily depends on the choice of $b$ in its $\sigma$-conjugacy class.
So we need to fix a specific representative to compare the $\J$-stratification on $X_\mu(b)$ (or $X_{\pc\mu}(b)$) to other stratification.
In loc. cit., it is pointed out that if $b$ is basic, then a reasonable choice is the unique length $0$ element in $B(G,\mu)$.
Also, for any $w\in \tW$, the $J_{\dot w}(F)$-stratification is independent of the choice of lift in $G(L)$.
See \cite[Lemma 2.5]{Gortz2}.

In the case where $G=\GL_n$ and $b=\tau^m$ with $m$ coprime to $n$, there is a group-theoretic way to describe the $\J$-stratification, which we will call the semi-module stratification.
Indeed, by \cite[Remark 3.1 \& Proposition 3.4]{CV}, the $\J$-stratification on $\cG r$ coincides with the stratification
$$G(L)/K=\bigsqcup_{\ld\in \Y} I \vp^\ld K/K.$$
So in this case, each $\J$-stratum of $X_\mu(b)$ (resp.\ $X_{\pc\mu}(b)$) coincides with $X_\mu^\ld(b)$ (resp.\ $X_{\pc\mu}^\ld(b)$) for some $\ld\in \Y$, where $X_\mu^\ld(b)=X_\mu(b)\cap I \vp^\ld K/K$ (resp.\ $X_{\pc\mu}^\ld(b)=X_{\pc\mu}(b)\cap I \vp^\ld K/K$).
Set $J_b(F)^0=J_b(F)\cap K=J_b(F)\cap I$.
Note that $\tau X_\mu^\ld(b)=X_\mu^{\tau\ld}(b)$ and  $J_b(F)/J_b(F)^0=\{\tau^k J_b(F)^0\mid k\in \Z\}$.
Thus $$J_b(F)X_\mu^\ld(b)=\bigsqcup_{k\in \Z} X_\mu^{\tau^k\ld}(b)\quad \text{and}\quad J_b(F)X_{\pc\mu}^\ld(b)=\bigsqcup_{k\in \Z} X_{\pc\mu}^{\tau^k\ld}(b).$$

See \S\ref{exsemi} for the precise definition of (extended) semi-modules. 
As we will explain in \S\ref{J-str semi-module}, the set $\{\ld\in\Y\mid X_\mu^\ld(b)\neq \emptyset\}$ can be regarded as semi-modules for $\mu$.
Let $w_{\max}$ be the longest element in $W_0$.
Then we have
$$\{\ld\in\Y\mid X_{-w_{\max}\mu}^\ld(b^{-1})\neq \emptyset\}=\{-w_{\max}\ld\mid \ld\in \cA_{\mu,b}\}.$$
Indeed it is easy to check that the image of $X_\mu^\ld(b)$ under the automorphism of $\cG r$ by $gK\mapsto w_{\max}{^tg}^{-1}K$ is $X_{-w_{\max}\mu}^{-w_{\max}\ld}(b^{-1})$.
This gives the description of ``dual" semi-modules for $\mu$.

\subsection{Length Positive Elements}
\label{LP}
We denote by $\delta^+$ the indicator function of the set of positive roots, i.e.,
$$\delta^+\colon \Phi\rightarrow \{0,1\},\quad \alpha \mapsto
\begin{cases}
1 & (\alpha\in \Phi_+) \\
0 & (\alpha\in \Phi_-).
\end{cases}$$
Note that any element $w\in \tW$ can be written in a unique way as $w=x\vp^\mu y$ with $\mu$ dominant, $x,y\in W_0$ such that $\vp^\mu y\in \SW$.
We have $p(w)=xy$ and $\ell(w)=\ell(x)+\la\mu, 2\rho\ra-\ell(y)$.
We define the set of {\it length positive} elements by $$\LP(w)=\{v\in W_0\mid \la v\alpha,y^{-1}\mu\ra+\delta^+(v\alpha)-\delta^+(xyv\alpha)\geq 0\  \text{for all $\alpha\in \Phi_+$}\}.$$
Then we always have $y^{-1}\in \LP(w)$.
Indeed, $y$ is uniquely determined by the condition that
$\la\alpha, \mu\ra\geq \delta^+(-y^{-1}\alpha)\ \text{for all $\alpha\in \Phi_+$}$.
Since $\delta^+(\alpha)+\delta^+(-\alpha)=1$, we have $$\la y^{-1}\alpha, y^{-1}\mu\ra+\delta^+(y^{-1}\alpha)-\delta^+(x\alpha)=\la \alpha,\mu\ra-\delta^+(-y^{-1}\alpha)+\delta^+(-x\alpha)\geq 0.$$
\begin{lemm}
\label{LPr}
For any $w=x\vp^\mu y\in \tW$ as above, we define $$\Phi_w\coloneqq\{\alpha\in \Phi_+\mid \la\alpha,\mu\ra-\delta^-(y^{-1}\alpha)+\delta^-(x\alpha)=0\}.$$
Here $\delta^-$ denotes the indicator function of the set of negative roots. Then we have
$$y\LP(w)=\{r^{-1}\in W_0\mid r(\Phi_+\setminus \Phi_w)\subset \Phi_+\ \text{or equivalently,}\ r^{-1}\Phi_+\subset \Phi_+\cup-\Phi_w\}.$$
\end{lemm}
\begin{proof}
Let $r\in W_0$ such that $r(\Phi_+\setminus \Phi_w)\subset \Phi_+$.
Let $\alpha\in \Phi_+$.
If $r^{-1}\alpha\in \Phi_+$, then we can check that $y^{-1}r^{-1}\in \LP(w)$ similarly as the case $r=1$ above.
If $r^{-1}\alpha\in \Phi_-$, then we must have $r^{-1}\alpha\in -\Phi_w$.
Since $\delta^-(-\alpha)=\delta^+(\alpha)$, it follows that
\begin{align*}
&\la y^{-1}r^{-1}\alpha,y^{-1}\mu\ra+\delta^+(y^{-1}r^{-1}\alpha)-\delta^+(xr^{-1}\alpha)\\
=&-(\la -r^{-1}\alpha, \mu\ra-\delta^-(-y^{-1}r^{-1}\alpha)+\delta^-(-xr^{-1}\alpha))=0.
\end{align*}
Thus $y^{-1}r^{-1}\in \LP(w)$.
This shows $\{r^{-1}\in W_0\mid r(\Phi_+\setminus \Phi_w)\subset \Phi_+\}\subseteq y\LP(w)$.

Let $v\in \LP(w)$ and let $\alpha\in \Phi_+$.
If $yv\alpha\in \Phi_-$, then
\begin{align*}
\la-yv\alpha,\mu\ra-\delta^-(-v\alpha)+\delta^-(-xyv\alpha)
=-(\la v\alpha, y^{-1}\mu\ra+\delta^+(v\alpha)-\delta^+(xyv\alpha))\le 0.
\end{align*}
On the other hand, by the characterization of $y$ above, we have
\begin{align*}
\la-yv\alpha,\mu\ra-\delta^-(-v\alpha)+\delta^-(-xyv\alpha)
=\la-yv\alpha,\mu\ra-\delta^+(v\alpha)+\delta^+(xyv\alpha)\geq 0.
\end{align*}
Thus $\la-yv\alpha,\mu\ra-\delta^-(-v\alpha)+\delta^-(-xyv\alpha)=0$ and hence $yv\alpha\in -\Phi_w$.
This shows $y\LP(w)\subseteq\{r^{-1}\in W_0\mid r(\Phi_+\setminus \Phi_w)\subset \Phi_+\}$.
The proof is finished.
\end{proof}

The notion of length positive elements is defined by Schremmer in \cite{Schremmer22}.
The description of $\LP(w)$ in Lemma \ref{LPr} is due to Lim \cite{Lim23}.

We say that the Dynkin diagram of G is $\sigma$-connected if it cannot be written as a union of two proper $\sigma$-stable subdiagrams that are not connected to each other. 
The following theorem is a refinement of the non-emptiness criterion in \cite{GHN2}, which is conjectured by Lim in \cite{Lim23} and proved by Schremmer in \cite[Proposition 5]{Schremmer23}.
\begin{theo}
\label{empty}
Assume that the Dynkin diagram of $G$ is $\sigma$-connected.
Let $b\in G(L)$ be a basic element with $\kappa(b)=\kappa(w)$.
Then $X_w(b)=\emptyset$ if and only if the following two conditions are satisfied:
\begin{enumerate}[(i)]
\item $|W_{\supp_\sigma(w)}|$ is not finite.
\item There exists $v\in \LP(w)$ such that $\supp_\sigma(\sigma^{-1}(v)^{-1}p(w)v)\subsetneq S$.
\end{enumerate}
\end{theo}

\begin{rema}
If $\kappa(b)\neq\kappa(w)$, then $X_w(b)=\emptyset$.
\end{rema}

\begin{rema}
Let $w\in \tW, w_0\in W_0$ and let $J\subseteq \Delta$ such that $J=\sigma(J)$.
Then we say that $w$ is a $(J,w_0,\sigma)$-{\it alcove element} if the following conditions are both satisfied:
\begin{enumerate}
\item $w_0^{-1}w\sigma(w_0)\in \tilde W_{J}\coloneqq \Y\rtimes W_{J}$, and
\item For any $\alpha\in w_0(\Phi_+\setminus \Phi_{J})$, $U_\alpha \cap {^wI}\subseteq U_\alpha\cap I$,
where $\Phi_J$ denotes the root system generated by $J$.
\end{enumerate}
In \cite[Proposition 5]{Schremmer23}, the condition (ii) in Theorem \ref{empty} is written as
\begin{enumerate}[(ii)']
\item There exist $J\subsetneq \Delta$ and $w_0\in W_0$ such that $w$ is a $(J,w_0,\sigma)$-alcove element.
\end{enumerate}
The equivalence of (ii) and (ii)' follows from \cite[Lemma 3.7 \& Lemma 3.9]{Lim23}.
\end{rema}

In the case $G=\GL_n$, there exists a length-preserving automorphism $\varsigma$ of $\tW$ defined as
$$w_0\vp^\ld\mapsto w_{\max}w_0w_{\max}^{-1}\vp^{-w_{\max}\ld},\quad w_0\in W_0,\ \ld\in \Y.$$
Note that $\varsigma(\tau^m)=\tau^{-m}$, $\varsigma(s_0)=s_0$ and $\varsigma(s_i)=s_{n-i}$ for $1\le i\le n-1$.
Let $w=x\vp^\mu y$ be as above.
For any $\alpha\in \Phi_+$ and $v\in \LP(w)$, we have
\begin{align*}
&\la \varsigma(v)(-w_{\max}\alpha), \varsigma(y^{-1})(-w_{\max}\mu)\ra+\delta^+(\varsigma(v)(-w_{\max}\alpha))-\delta^+(\varsigma(xy)\varsigma(v)(-w_{\max}\alpha))\\
&=\la v\alpha,y^{-1}\mu\ra+\delta^+(v\alpha)-\delta^+(xyv\alpha)\geq 0.
\end{align*}
Thus $\LP(\varsigma(w))=\varsigma(\LP(w))=w_{\max}\LP(w)w_{\max}^{-1}$.
In particular, there exists $v\in \LP(w)$ such that $v^{-1}p(w)v$ is a Coxeter element if and only if the same is true for $\varsigma(w)$ and $\LP(\varsigma(w))$.

\section{Semi-Modules}
\label{semi-modules}
From now and until the end of this paper, we set $G=\GL_n$ and $b=\tau^m$.
For $\mu\in \Y_+$, let $\mu(i)$ denotes the $i$-th entry of $\mu$.
Then $[\tau^m]\in B(G,\mu)$ if and only if $m=\mu(1)+\cdots+\mu(n)$.
We assume this from now.
Also, without loss of generality, we may and will assume that $\mu(n)=0$.
Recall that $w_{\max}$ is the longest element in $W_0$.

\subsection{Extended Semi-Modules}
\label{exsemi}
Here we recall the definition of extended semi-modules in a combinatorial way from \cite{Viehmann}.
Note that although we choose the subgroup of upper triangular matrices $B$ as Borel subgroup in this paper, the fixed Borel subgroup in \cite{Viehmann} is the subgroup of {\it lower} triangular matrices. 
\begin{defi}
\label{semi-module}
A {\it semi-module} for $m, n$ is a subset $A\subset \Z$ that is bounded below and satisfies $m+A\subset A$ and $n+A\subset A$.
Set $\bar{A}=A\setminus (n+A)$.
The semi-module $A$ is called normalized if $\sum_{a\in \bar{A}}a=\frac{n(n-1)}{2}$.
\end{defi}

For a semi-module $A$, there exists a unique $\mu'\in \N^n$ satisfying the following condition:
Let $a_0=\min \bar A$ and let inductively $a_i=a_{i-1}+m-\mu'(i)n$ for $i=1,\ldots, n$.
Then $a_0=a_n$ and $\{a_0,a_1,\ldots,a_{n-1}\}=\bar A$.
We call $\mu'$ the {\it type} of $A$.

\begin{lemm}
\label{type}
There is a bijection between the set of normalized semi-modules for $m,n$ and the set of possible types $\mu'\in \N^n$ with $\nu_b\pc w_{\max}\mu'$.
\end{lemm}
\begin{proof}
This is \cite[Lemma 3.3]{Viehmann}.
\end{proof}

\begin{defi}
\label{exsemidef}
An {\it extended semi-module} $(A,\varphi)$ for $\mu\in \Y_+$ is a normalized semi-module $A$ for $m,n$ together with a function $\vph\colon \Z\rightarrow \N\cup\{-\infty\}$ satisfying the following properties:
\begin{enumerate}[(1)]
\item $\vph(a)=-\infty$ if and only if $a\notin A$.
\item $\vph(a+n)\geq \vph(a)+1$ for all $a\in \Z$.
\item $\vph(a)\le \max\{k\mid a+m-kn\in A\}$ for all $a\in A$.
If $b\in A$ for all $b\geq a$, then the two sides are equal.
\item There is a decomposition of $A$ into disjoint union of sequences $a_j^1,\ldots, a_j^n$ with $j\in \N$ and the following properties:
\begin{enumerate}
\item $\vph(a_{j+1}^l)=\vph(a_j^l)+1$.
\item If $\vph(a_j^l+n)=\vph(a_j^l)+1$, then $a_{j+1}^l=a_j^l+n$. Otherwise $a_{j+1}^l>a_j^l+n$.
\item The $n$-tuple $(\vph(a_0^l))$ is a permutation of $\mu$.
\end{enumerate}
\end{enumerate}
An extended semi-module such that the equality holds in (3) for all $a\in A$ is called {\it cyclic}.
\end{defi}

For any $\ld\in \Y$, we denote by $\ld_{\dom}$ the dominant conjugate of $\ld$.
Let $\mu'$ be the type of a semi-module for $m,n$.
Let $\vph$ be a function such that $(1)$ and the equation in $(3)$ hold.
Then it is easy to check that $(A,\vph)$ is a cyclic semi-module for $\mu'_{\dom}$.
In general, the following lemma holds.
\begin{lemm}
\label{cyclic}
Let $(A,\vph)$ be an extended semi-module for $\mu$ and let $\mu'$ be the type of $A$.
Then $\mu'_{\dom}\pc \mu$ and $(A,\vph)$ is cyclic if and only if $\mu'\in W_0\mu$.
Moreover, if $\mu$ is minuscule, then all extended semi-modules for $\mu$ are cyclic.
\end{lemm}
\begin{proof}
See \cite[Lemma 3.6 \& Corollary 3.7]{Viehmann}.
See also \cite[Lemma 5.9]{Hamacher}.
\end{proof}

Let $e_0,\ldots, e_{n-1}$ be the standard basis of $L^n$.
Then the lattice $\cO^n$ is generated by $e_0,\ldots, e_{n-1}$.
For $i\in \Z$, we define $e_i$ by $e_{i+n}=\vp e_i$.
Note that we have $\tau e_i=e_{i+1}$ for any $i$.
In the sequel, we identify $\cG r$ and $\{M\subset L^n\ \text{lattice}\}$ by $gK\mapsto g\cO^n$.

Let $X_\mu(b)^0$ be a $\aFq$-subscheme of $X_\mu(b)$ defined as $X_\mu(b)^0=\{gK\in X_\mu(b)\mid \kappa(g)=0\}$.
We associate to $M\in X_\mu(b)^0$ an extended semi-module for $\mu$.
Let $v\in L^n$.
Then we can write $v=\sum_{i\in \Z}\alpha_ie_i$ with $\alpha_i\in \aFq$ and $\alpha_i=0$ for sufficiently small $i$.
Let $$\cI\colon L^n\setminus \{0\}\rightarrow \Z,\quad v\mapsto \min\{i\mid \alpha_i\neq 0\}.$$
For $M\in \cG r$, we define the set 
$$A(M)=\{\cI(v)\mid v\in M\setminus \{0\}\}.$$
It is easy to check that if $M\in X_\mu(b)^0$, then $A(M)$ is a normalized semi-module for $m,n$.
We also define $\vph(M)\colon \Z\rightarrow \N\cup \{-\infty\}$ by
$$a\mapsto
\begin{cases}
\max\{k\mid \exists v\in M \text{\ with\ }\cI(v)=a, \vp^{-k}b\sigma(v)\in M\} & (a\in A(M))\\
-\infty & (a\notin A(M)).
\end{cases}$$
\begin{lemm}
Let $M\in X_\mu(b)^0$.
Then $(A(M), \vph(M))$ is an extended semi-module for $\mu$.
\end{lemm}
\begin{proof}
See \cite[Lemma 4.1]{Viehmann}.
\end{proof}

For an extended semi-module $(A,\vph)$ for $\mu$, let
$$S_{A,\vph}=\{M\mid A(M)=A,\vph(M)=\vph\}\subset \cG r.$$

\begin{lemm}
The set $S_{A,\vph}$ is a locally closed subscheme of $X_\mu(b)^0$.
\end{lemm}
\begin{proof}
See \cite[Lemma 4.2]{Viehmann}.
\end{proof}

Let $\A_{\mu,b}$ be the set of extended semi-modules for $\mu$.
Set $\A_{\mu,b}^\tp= \{(A,\varphi)\in \A_{\mu,b}\mid \dim S_{A,\varphi}=\dim X_\mu(b)\}$.
By Proposition \ref{UAphi} below, $J_b(F)\backslash \Irr X_\mu(b)$ is parametrized by $\A_{\mu,b}^\tp$.
In the sequel, we also use the symbol $\A$ to denote the affine space as usual.
We hope our notation will not cause confusions.

For an extended semi-module $(A,\vph)$ for $\mu$, let
$$\cV(A,\vph)=\{(a,c)\in A\times A\mid c>a, \vph(a)>\vph(c)>\vph(a-n)\}.$$
\begin{prop}
\label{UAphi}
Let $(A,\vph)$ be an extended semi-module for $\mu$.
There exists a non-empty open subscheme $U_{A,\vph}\subseteq \A^{|\cV(A,\vph)|}$ and a morphism $U_{A,\vph}\rightarrow S_{A,\vph}$ which is bijective on $\aFq$-valued points.
In particular, $S_{A,\vph}$ is irreducible and of dimension $|\cV(A,\vph)|$.
Moreover if $(A,\vph)$ is a cyclic extended semi-module, then $U_{A,\vph}=\A^{|\cV(A,\vph)|}$.
\end{prop}
\begin{proof}
See \cite[Theorem 4.3]{Viehmann}.
\end{proof}

Here we briefly describe $U_{A,\vph}$ and the map $U_{A,\vph}\rightarrow S_{A,\vph}$.
For any $x\in \aFq^{|\cV(A,\vph)|}=\A^{|\cV(A,\vph)|}$, we denote the coordinates of $x$ by $x_{a,c}$.
We associate to every $x$ a set of elements $\{v(a)\in L^n\mid a\in A\}$ which satisfies the following equations.

If $a=\max\bar A$, then
$$v(a)=e_a+\sum_{(a,c)\in \cV(A,\vph)}[x_{a,c}] v(c).$$
For any other element $a\in \bar A$, we want
$$v(a)=v'+\sum_{(a,c)\in \cV(A,\vph)}[x_{a,c}] v(c),$$
where $v'=\vp^{-\vph(a')}b\sigma(v(a'))$ for $a'$ being minimal satisfying $a'+m-\vph(a')n=a$.
For $a\in n+A$, we want
$$v(a)=\vp v(a-n)+\sum_{(a,c)\in \cV(A,\vph)}[x_{a,c}] v(c).$$
Here $[x_{a,c}]$ denotes the Teichm\"{u}ller lift of $x_{a,c}$ if $\ch F=0$ and $[x_{a,c}]=x_{a,c}$ if $\ch F>0$.
The set $\{v(a)\in L^n\mid a\in A\}$ is uniquely determined by the equations above.
Hence the map $\A^{|\cV(A,\vph)|}\rightarrow \cG r, x\mapsto \la v(a)\ra_{a\in A}$ is well-defined.
By applying $\sigma$ on the above equations for $x$, we can easily check that this map is compatible with the action of $\sigma$, i.e., $\sigma(x)\coloneqq(x_{a,c}^q)$ maps to $\sigma\la v(a)\ra_{a\in A}$.
Let $U_{A,\vph}$ be the preimage of $S_{A,\vph}$ under this map.
Then $S_{A,\vph}$ and hence $U_{A,\vph}$ are stable under $\sigma$ (because $\sigma(b)=b$).
In particular, we have $|S_{A,\vph}^\sigma|=|U_{A,\vph}^\sigma|$.
So if $(A,\vph)$ is cyclic, then  $|S_{A,\vph}^\sigma|=q^{|\cV(A,\vph)|}$.
Although not needed in this paper, it is also worth mentioning that if $(A,\vph)$ is non-cyclic, then $S_{A,\vph}$ is never universally homeomorphic to an affine space.

\begin{prop}
\label{affine}
If $(A,\vph)$ is non-cyclic, then $|S_{A,\vph}^\sigma|<q^{|\cV(A,\vph)|}$.
In particular, $S_{A,\vph}$ is never universally homeomorphic to an affine space.
\end{prop}
\begin{proof}
Let $x\in \A^{|\cV(A,\vph)|}$.
Note that if $x_{a,c}=0$ for all $(a,c)\in \cV(A,\vph)$, then $v(a)=e_a$ for all $a\in A$.
Set $M=\la e_a\ra_{a\in A}$.
Then it is easy to check that $(A(M),\vph(M))$ is a cyclic semi-module for the dominant conjugate of the type of $A(M)$.
So if $(A,\vph)$ is not cyclic, then $M\notin S_{A,\vph}$ and hence $|S_{A,\vph}^\sigma|=|U_{A,\vph}^\sigma|<q^{|\cV(A,\vph)|}$.
The last statement follows from \cite[Proposition 4.1.12 \& Proposition 8.1.11 (ii)]{DM}.
\end{proof}

\subsection{The Stratification by Extended Semi-Modules}
\label{J-str semi-module}
For any $\ld\in \Y$, set $A^\ld=\{(i-1)+\ld(i)n+\N n \mid 1\le i\le n\}.$
It is easy to check that for a lattice $M\in I\vp^\ld K/K$, we have $A(M)=A^\ld$.
Thus we have the following lemma, which relates the semi-module stratification to the stratification by extended semi-modules.
\begin{lemm}
\label{semi-module str}
Let $\ld\in \Y$ with $\ld(1)+\cdots+\ld(n)=0$.
Then $X_\mu^{\ld}(b)\neq \emptyset$ if and only if there exists an extended semi-module $(A^\ld,\vph)$ for $\mu$.
If this is the case, we have
$$X_\mu^{\ld}(b)=\bigsqcup_{\vph}S_{A^\ld,\vph},$$
where $\vph$ runs over all the functions $\Z\rightarrow \N\cup\{-\infty\}$ such that the pair of $A^\ld$ and the function is an extended semi-module for $\mu$.
\end{lemm}

For $\ld\in \Y$ with $X_\mu^\ld(b)\neq \emptyset$, let $1\le i_0\le n$ such that $(i_0-1)+\ld(i_0)n=\min \overline {A^\ld}$.
Then 
\begin{align*}
&(i_0-1)+\ld(i_0)n+m-(\ld(i_0)+\ld_{b,\dom}(c^m(i_0))-\ld(c^m(i_0)))n\\
=& c^m(i_0)-1+\ld(c^m(i_0))n\in \bar A,
\end{align*}
where $c=s_1\cdots s_{n-1}$.
Repeating the same argument, we can check that the type of $A^\ld$ is a conjugate of $b\ld-\ld=c^m\ld+\ld_{b,\dom}-\ld$.
By Lemma \ref{cyclic}, an extended semi-module $(A^\ld,\vph)$ for $\mu$ is cyclic if and only if $b\ld-\ld\in W_0\mu$.

\begin{coro}
\label{never refinement}
Let $\mu\in \Y_+$.
If there exists a non-cyclic semi-module for $\mu$, then the semi-module stratification of $X_{\pc\mu}(b)$  is not a refinement of the Ekedahl-Oort stratification.
\end{coro}
\begin{proof}
Let $(A^\ld,\vph)$ be a non-cyclic semi-module for $\mu$.
Then we have $(b\ld-\ld)_{\dom}\prec \mu$ by Lemma \ref{cyclic}.
On the other hand, there always exists a cyclic semi-module $(A^\ld,\vph')$ for $(b\ld-\ld)_{\dom}$.
By Lemma \ref{semi-module str}, $X_{\pc\mu}^\ld(b)$ intersects both $X_\mu(b)$ and $X_{(b\ld-\ld)_{\dom}}(b)$.
This implies that $X_{\pc\mu}^\ld(b)$ is not contained in any set of the form $\pi(X_w(b))$ with $w\in \tW$, which finishes the proof.
\end{proof}

For $\mu=(\mu(1),\ldots, \mu(n-1),0)\in X_*(T)_+$, set $\mu^*=(\mu(1), \mu(1)-\mu(n-1),\ldots,\mu(1)-\mu(2),0)$ and $b^*=\tau^{n\mu(1)-m}$.
If $(A^\ld,\vph)$ is an extended semi-module for $\mu$, then there exists $\vph'\colon \Z\rightarrow \N\cup\{-\infty\}$ such that $(A^{-w_{\max}\ld},\vph')$ is an extended semi-module for $\mu^*$ (see \S\ref{J-str}).
Clearly, $b\ld-\ld\in W_0\mu$ if and only if $b^*(-w_{\max}\ld)+w_{\max}\ld\in -W_0\mu^*$.
Thus we have the following lemma.
\begin{lemm}
\label{dual cyclic}
There exists a non-cyclic extended semi-module for $\mu$ if and only if the same is true for $\mu^*$.
\end{lemm}

\subsection{The Minuscule Case}
\label{minuscule case}
In this subsection, we treat the minuscule case.
Consider $G^d$ with a Frobenius automorphism given by
$$(g_1,g_2,\ldots, g_d)\mapsto (g_2,\ldots, g_d, \sigma(g_1)).$$
For $\mu_\bl=(\mu_1, \ldots, \mu_d)\in \Y^d_+$ and $b_\bl=(1,\ldots, 1, b)\in G^d(L)$ with $b\in G(L)$, we  define $X_{\mu_\bl}(b_\bl)\subset \cG r^d=G^d(L)/K^d$ as
$$X_{\mu_\bl}(b_\bl)=\{x_\bl K^d\in \cG r^d\mid x_\bl^{-1}b_\bl\sigma(x_\bl)\in K^d\vp^{\mu_\bl}K^d\}.$$

Let us denote by $\Irr X_{\mu_\bl}(b_\bl)$ the set of irreducible components of $X_{\mu_\bl}(b_\bl)$.
Through the identification $J_b(F)\cong J_{b_\bl}(F)$ given by $g\mapsto (g,\ldots,g)$, this set is equipped with an action of $J_b(F)$.

For minuscule $\mu_\bl\in \Y^d_+$ and $b_\bl=(1,\ldots,1,b)\in G^d(L)$,
we define
$$\cA_{\mu_\bl, b_\bl}^{\tp}\coloneqq \{\ld_\bl\in \Y^d\mid \dim X_{\mu_\bl}^{\ld_\bl}(b_\bl)=\dim X_{\mu_\bl}(b_\bl)\}.$$
Here $X_{\mu_\bl}^{\ld_\bl}(b_\bl)$ denotes $X_{\mu_\bl}(b_\bl)\cap It^{\ld_\bl} K/K$.
For $\ld_\bl, \ld'_\bl\in \cA_{\mu_\bl, b_\bl}^{\tp}$, we write $\ld_\bl\sim \ld'_\bl$ if $\ld_\bl=\tau^k \ld'_\bl=(\tau^k \ld'_1,\ldots, \tau^k \ld'_d)$ for some $k\in \Z$.
Let $\A_{\mu_\bl, b_\bl}^{\tp}$ denote the set of equivalence classes with respect to $\sim$, and let $[\ld_\bl]\in \A_{\mu_\bl, b_\bl}^{\tp}$ denote the equivalence class represented by $\ld_\bl\in \cA_{\mu_\bl, b_\bl}^{\tp}$.
Then $J_{b}(F)\backslash \Irr X_{\mu_\bl}(b_\bl)$ is parametrized by $\A_{\mu_\bl, b_\bl}^{\tp}$ as follows.
\begin{prop}
\label{minuscule bijection}
Assume that $\mu_\bl\in \Y_+^d$ is minuscule.
Then the map $\ld_\bl\mapsto \overline{X_{\mu_\bl}^{\ld_\bl}(b_\bl)}$ induces a bijection
$$\A_{\mu_\bl,b_\bl}^\tp\cong J_{b}(F)\backslash \Irr X_{\mu_\bl}(b_\bl).$$
\end{prop}
\begin{proof}
See \cite[Proposition 1.6]{HV2}.
Note that we have $\mathrm{Stab}_{J_b(F)} (X_{\mu_\bl}^{\ld_\bl}(b_\bl))=J_b(F)^0$.
\end{proof}

We also define
$$\cA_{\mu_\bl, b_\bl}^j\coloneqq \{\ld_\bl\in \Y^d\mid  \dim X_{\mu_\bl}^{\ld_\bl}(b_\bl)=j\}$$
for $1\le j\le \dim X_{\mu_\bl}(b_\bl)$.
We can similarly consider the equivalence relation $\sim$ as above.
If $d=1$, then $\A_{\mu,b}^j\coloneqq\cA_{\mu, b}^j/\sim$ can be identified with (extended) semi-modules for $\mu$ whose corresponding stratum has dimension $i$, see Lemma \ref{cyclic} and Lemma \ref{semi-module str}.

\begin{prop}
\label{semi-module number}
Set $\mu=\omega_i$.
Then we always have $|\A_{\mu,b}^{\tp}|=|\A_{\mu,b}^{0}|=1$.
If $i=2,n-2$, then $|\A_{\mu,b}^j|=1$ for all $0\le j\le \dim X_{\mu}(b)$.
If $i=3,n-3$, then $|\A_{\mu,b}^{\dim X_{\mu}(b)-1}|=2$.
\end{prop}
\begin{proof}
We can easily check the equalities in the proposition using \cite[Theorem 4.16]{HV2} (cf.\ \cite[Remark 6.16]{JO}), which gives a combinatorial way of computing $|\A_{\mu,b}^j|$.
In fact, all of the assertions except the last assertion follow from \cite[Proposition 5.5]{Viehmann2}.
\end{proof}

\begin{exam}
We always have $\A_{\omega_i,b}^{0}=\{[0]\}$.
\end{exam}

\section{Crystal Bases}
\label{crystal bases}
Keep the notations and assumptions in \S\ref{semi-modules}.
\subsection{Crystals and Young Tableaux}
In this subsection, we first recall the definition of $\hG$-crystals from \cite[Definition 3.3.1]{XZ}.
\begin{defi}
\label{crystaldefi}
A (normal) {\it $\hG$-crystal} is a finite set $\B$, equipped with a weight map $\wt\colon \B\rightarrow \Y$, and operators $\te_\alpha, \tf_\alpha\colon \B\rightarrow \B\cup \{0\}$ for each $\alpha\in \Delta$, such that
\begin{enumerate}[(i)]
\item for every $\bb\in \B$, either $\te_\alpha\bb=0$ or $\wt(\te_\alpha \bb)=\wt(\bb)+\alpha^\vee$, and either $\tf_\alpha\bb=0$ or $\wt(\tf_\alpha \bb)=\wt(\bb)-\alpha^\vee$,
\item for all $\bb, \bb'\in \B$ one has $\bb'=\te_\alpha \bb$ if and only if $\bb=\tf_\alpha \bb'$, and
\item if $\ve_\alpha, \phi_\alpha\colon \B\rightarrow \Z,\ \alpha\in \Delta$ are the maps defined by
\begin{align*}
\ve_\alpha(\bb)=\max\{k\mid \te_\alpha^k\bb\neq 0\}\ \  \text{and}\ \ \phi_\alpha(\bb)=\max\{k\mid \tf_\alpha^k\bb\neq 0\},
\end{align*}
then we require $\phi_\alpha(\bb)-\ve_\alpha(\bb)=\la \alpha, \wt(\bb)\ra$.
\end{enumerate}
For a $\hG$-crystal $\B$, let $\B^*=\{\bb^*\mid \bb\in \B\}$ be the dual $\hG$-crystal.
Setting $0^*=0$, the maps are given by 
$$\wt(\bb^*)=-\wt(\bb),\quad \te_\alpha(\bb^*)=(f_\alpha \bb)^*,\quad \text{and}\quad \tf_\alpha(\bb^*)=(\te_\alpha \bb)^*.$$
\end{defi}

For $\ld\in \Y$, we denote by $\B(\ld)$ the set of elements with weight $\ld$ for $\hG$, called the {\it weight space} with weight $\ld$ for $\hG$.
Let $\B_1$ and $\B_2$ be the two $\hG$-crystals.
A morphism $\B_1\rightarrow \B_2$ is a map of underlying sets compatible with $\wt,\te_\alpha$ and $\tf_\alpha$.

In the sequel, we write $\te_i$ and $\tf_i$ (resp.\ $\ve_i$ and $\phi_i$) instead of $\te_{\chi_{i,i+1}}$ and $\tf_{\chi_{i,i+1}}$ (resp.\ $\ve_{\chi_{i,i+1}}$ and $\phi_{\chi_{i,i+1}}$) for simplicity.

\begin{exam}
\label{highest}
Let $\B_\mu$ be the crystal basis of the irreducible $\hG$-module of highest weight $\mu\in \Y_+$.
Then $\B_\mu$ is a crystal.
We call $\B_\mu$ a {\it highest weight crystal} of highest weight $\mu$ (cf.\ \cite[Definition 3.3.1 (3)]{XZ}).
There exists a unique element $\bb_\mu\in \B_\mu$ satisfying $\te_\alpha\bb_\mu=0$ for all $\alpha$, $\wt(\bb_\mu)=\mu$, and $\B_\mu$ is generated from $\bb_\mu$ by operators $\tf_\alpha$.
So $w_{\max}\bb_{\mu}$ has the {\it lowest} weight $w_{\max}\mu$.
It is well-known that the dual of $\B_\mu$ is isomorphic to $\B_{-w_{\max}\mu}$ (see for example \cite[Lemma 3.5.2]{HK}).
\end{exam}

We give a realization of $\B_\mu$ by Young tableaux.
This allows us to treat it in a combinatorial way.

\begin{defi}
A {\it Young diagram} is a collection of boxes arranged in left-justified rows with a weakly decreasing number of boxes in each row.
For a dominant cocharacter $\mu\in \Y_+$, we denote by $Y_\mu$ the Young diagram having $\mu(i)$ boxes in the $i$th row.
A {\it skew Young diagram} is a diagram obtained by removing a smaller Young diagram from a larger one that contains it.
For dominant cocharacters $\mu,\nu\in \Y_+$ with $\nu(i)\le \mu(i)$, we denote by $Y_{\mu/\nu}$ the skew Young diagram obtained by removing $Y_\nu$ from $Y_\mu$. 
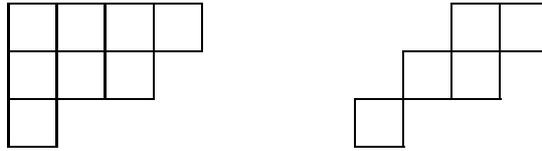
$$
\ytableausetup{nosmalltableaux}
\begin{ytableau}
   \null & \null & \null & \null \\
   \null & \null   & \null \\
   \null
 \end{ytableau}
\hspace{2cm}
\ytableausetup{nosmalltableaux}
\begin{ytableau}
   \none & \none & \null & \null \\
   \none &   \null   & \null \\
   \null
\end{ytableau}
$$
\end{defi}

\begin{defi}
A {\it tableau} is a (skew) Young diagram filled with numbers, one for each box.
A {\it semi-standard tableau} is a tableau obtained from a (skew) Young diagram by filling the boxes with the numbers $1,2,\ldots, n$ subject to the conditions
\begin{enumerate}[(i)]
\item the entries in each row are weakly increasing from left to right,
\item the entries in each column are strictly increasing from top to bottom.
\end{enumerate}
$$
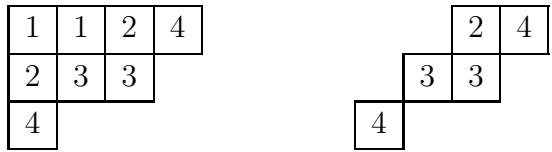
\begin{ytableau}
   1 & 1 & 2 & 4 \\
   2 & 3 & 3 \\
   4
  \end{ytableau}
\hspace{2cm}
\ytableausetup{nosmalltableaux}
\begin{ytableau}
   \none & \none & 2 & 4 \\
   \none &   3   & 3 \\
   4
  \end{ytableau}
$$
\end{defi}

Let $K_{\mu/\nu}(\ld)$ be the number of all semi-standard tableaux $\bb$ of shape $Y_{\mu/\nu}$ such that the number of $\young(i)$ appearing in $\bb$ is $\ld(i)$ for $1\le i\le n$.
This is sometimes called the {\it Kostka number}.
In \S\ref{NCtableaux}, we need the following well-known result.
\begin{prop}
\label{Kostka}
Let $\ld,\ld'\in \Y_+$.
If $\ld\pc \ld'$, then $K_{\mu/\nu}(\ld')\le K_{\mu/\nu}(\ld)$.
In particular, $K_{\mu/\nu}(\ld')\neq 0$ implies $K_{\mu/\nu}(\ld)\neq 0$.
\end{prop}
\begin{proof}
See \cite[Proposition 1.2]{Fayers} and the remark right after the proposition.
\end{proof}

We denote by $\cB(Y)$ the set of all semi-standard tableaux of shape $Y$.

\begin{theo} 
\label{Young theorem}
Let $\mu=(\mu(1),\ldots,\mu(n))\in \Y_+\setminus\{0\}$ with $\mu(n)=0$.
Then $\cB(Y_\mu)$ has a crystal structure.
Moreover, the crystal $\cB(Y_\mu)$ is isomorphic to $\B_\mu$.
\end{theo}
\begin{proof}
This is \cite[Theorem 7.3.6 \& Theorem 7.4.1]{HK}.
\end{proof}

In the sequel, we identify $\B_\mu$ and $\cB(Y)$ by Theorem \ref{Young theorem}.
For a semi-standard tableau $\bb\in \B_\mu$, let $k_i$ denote the number of $i$'s appearing in $\bb$.
Then the weight map $\wt$ on $\B_\mu$ is given by $\wt(\bb)=(k_1,\ldots, k_n)$.
The following result is an explicit description of the actions of $\te_i$ and $\tf_i$ on $\B_\mu$.

\begin{theo}
\label{ef}
The actions of $\te_i$ and $\tf_i$ on $\bb\in \B_\mu$ can be computed by following the steps below:
\begin{enumerate}[(i)]
\item In the Far-Eastern reading $\bb_1\otimes\cdots\otimes \bb_N$ of $\bb$,
we identify \framebox[1.2em]{\rule{0pt}{1.6ex}$i$} (resp.\ \fbox{$i+1$}) by $+$ (resp.\ $-$) and neglect other boxes. 
\item Let $u_i(\bb)=u^1u^2\cdots u^l\ (u^j\in \{\pm\})$ be the sequence obtained by (i).
If there is ``$+-$'' in $u(\bb)$, then we neglect such a pair.
We continue this procedure as far as we can.
\item Let $u_i(\bb)_{\red}=-\cdots-+\cdots+$ be the sequence obtained by (ii).
Then $\te_i$ changes the rightmost $-$ in $u(\bb)_{\red}$ to $+$, and $\tf_i$ changes the leftmost $+$ in $u(\bb)_{\red}$ to $-$.
If there is no such $-$ (resp.\ $+$), then $\te_i\bb=0$ (resp.\ $\tf_i \bb=0$).
\end{enumerate}
Moreover, $\ve_i(\bb)$ (resp.\ $\phi_i(\bb)$) is equal to the number of $-$ (resp.\ $+$) in $u(\bb)_{\red}$.
\end{theo}
\begin{proof}
The first statement is \cite[Theorem 3.4.2]{KN}.
The second statement follows immediately from this.
\end{proof}

Next we recall the Weyl group action on crystals.
Let $\B$ be a $\hG$-crystal.
For any $1\le i\le n-1$ and $\bb\in \B$, we set
\begin{align*}
s_i\bb=
\begin{cases}
\tf_i^{\la \chi_{i, i+1}, \wt(\bb)\ra}\bb & \text{if $\la \chi_{i, i+1}, \wt(\bb)\ra \geq 0$}\\
\te_i^{-\la \chi_{i, i+1}, \wt(\bb)\ra}\bb & \text{if $\la \chi_{i, i+1}, \wt(\bb)\ra \le 0$}.
\end{cases}
\end{align*}
Then we have the obvious relation
$$\wt(s_i\bb)=s_i(\wt(\bb)).$$
By \cite[Theorem 7.2.2]{Kashiwara}, this extends to the action of the Weyl group $W_0$ on $\B$, which is compatible with the action on $\Y$.

\begin{lemm}
\label{Weyl action}
Let $w, w'\in W_0$ and $\bb\in \B$.
If $w(\wt(\bb))=w'(\wt(\bb))$, then $w\bb=w'\bb$.
\end{lemm}
\begin{proof}
This is \cite[Lemma 3.10]{Shimada3}.
\end{proof}
Let $\bb\in \B(\ld)$.
If $\ld'$ is a conjugate of $\ld$, i.e., there exists $w\in W_0$ such that $\ld'=w\ld$, then we call $w\bb$ the conjugate of $\bb$ with weight $\ld'$.
By Lemma \ref{Weyl action}, this does not depend on the choice of $w$.

Finally we consider the minuscule case.
If $\mu\in \Y_+$ is minuscule, then $\wt\colon\B_\mu\rightarrow \Y$ gives an identification between $\B_\mu$ and the set of cocharacters which are conjugate to $\mu$.
Suppose $\mu_\bl=(\mu_1,\ldots,\mu_d)\in \Y_+^d$ is minuscule.
We can also identify $\B_{\mu_\bl}^{\hG^d}:=\B_{\mu_1}\times\cdots\times \B_{\mu_d}$ with the set of cocharacters in $\Y^d$ which are conjugate to $\mu_\bl$.

For $1\le k<n$, let $\omega_k$ be the cocharacter of the form $(1,\ldots,1,0,\ldots,0)$ in which $1$ is repeated $k$ times.
Assume that each $\mu_i$ is equal to $\omega_{k_i}$ for some $1\le k_i<n$ and $i<j$ if and only if $k_i\le k_j$.
In the rest of paper, we call such $\mu_\bl$ {\it Far-Eastern}.
If $\mu_\bl$ is Far-Eastern, then $|\mu_\bl|:=\mu_1+\cdots+\mu_d$ is dominant and its last entry is $0$.
Let $\FE\colon \B_{|\mu_\bl|}\rightarrow \B_{\mu_\bl}^{\hG^d}$ be a map defined by decomposing $\bb\in \B_\mu$ into its columns from right to left.
We call $\FE$ the Far-Eastern reading.

\subsection{Construction of Extended Semi-Modules}
In this subsection, we recall from \cite[\S 4.2]{Shimada3} the way of constructing extended semi-modules.
See \cite[\S 4.3]{Shimada3} for some examples of computation.
Let $\mu_\bl\in \Y^d_+$ be a Far-Eastern cocharacter.
Set $\mu=|\mu_\bl|$.

Let $\ld_b$ denote the cocharacter whose $i$-th entry is $\lfloor \frac{im}{n}\rfloor-\lfloor \frac{(i-1)m}{n}\rfloor$.
Set $\ld_b^{\op}=w_{\max}\ld_b$.
For any $\bb\in \B_\mu(\ld_b)$, we denote by $\bb^{\op}$ the conjugate of $\bb$ with weight $\ld_b^{\op}$.
Let $1\le m_0<n$ be the residue of $m$ modulo $n$.
Note that each entry of $\ld_b$ is $\mn$ or $\mn+1$ and $\ld_b(i)=\ld_b(n+1-i)$ for any $2\le i\le n-1$.
Let $i_0=1<i_1<i_2<\cdots<i_{m_0}=n$ be the integers such that $\ld_b(i_1)=\ld_b(i_2)=\cdots =\ld_b(i_{m_0})=\mn+1$.
Then 
$$\text{$\ld_b^{\op}=w_{\max}'\ld_b$,\quad where $w_{\max}'=(s_{i_{m_0-1}}\cdots s_{n-1})\cdots (s_{i_1}\cdots s_{i_2-1})(s_1\cdots s_{i_1-1})$.}$$
Here $\ld_b(i)=\mn$ (resp.\ $\ld_b(i+1)=\mn$) if and only if $s_{i-1}s_i\le w_{\max}'$ (resp.\ $s_is_{i+1}\le w_{\max}'$).
By Lemma \ref{Weyl action}, it follows that $\bb^{\op}$ can be computed by the action of the Coxeter element $w_{\max}'$.
In this computation, each $s_i$ acts as the action of $\te_i$ because $\mn-(\mn+1)=-1$.
Therefore, if we write $$\FE(\bb)=(\bb_1,\ldots,\bb_d)$$
then there exists $(w_1, \ldots, w_{d})\in W_0^d$ such that
$$\FE(\bb^{\op})=(w_1\bb_1, \ldots, w_d\bb_d)$$
and each simple reflection appears exactly once in some $\supp(w_j)$.
\begin{lemm}
\label{w}
The tuple $(w_1,\ldots,w_d)\in W_0^d$ as above is uniquely determined by $\bb$.
In particular, $w(\bb)\coloneqq w_1^{-1}\cdots w_d^{-1}$ is a Coxeter element uniquely determined by $\bb$.
\end{lemm}
\begin{proof}
This is \cite[Lemma 4.3]{Shimada3}.
\end{proof}
Set $w(\bb)=w_1^{-1}\cdots w_d^{-1}$ and $\Upsilon(\bb)=\{\up\in W_0\mid \up^{-1}c^m \up=w(\bb)\}$, where $c=s_1s_2\cdots s_{n-1}$.
Clearly $|\Upsilon(\bb)|=n$.

For any $\bb'\in \B_{\mu}$, set
$$\xi(\bb')=(\ve_1(\bb')+\cdots+\ve_{n-1}(\bb'),\ve_2(\bb')+\cdots+\ve_{n-1}(\bb'),\ldots,\ve_{n-1}(\bb'),0).$$
Let $\ld_b^-$ be the anti-dominant conjugate of $\ld_b$, and let $\bb^-$ be the conjugate of $\bb$ with weight $\ld_b^-$.
For any $\bb\in \B_{\mu}(\ld_b)$ and $\up\in \Upsilon(\bb)$, we define $\xi_{\bl}(\bb, \up)\in \Y^d$ by
$$\xi_{j}(\bb, \up)=\up \xi(\up^{-1}\bb^-)+\sum_{1\le j'< j}\up w_1^{-1}\cdots w_{j'-1}^{-1}\wt(\bb_{j'})\quad (1\le j\le d).$$

Let $C\in \Irr X_\mu(b)^0$.
By Proposition \ref{UAphi}, $C=\overline{S_{A,\vph}}$ for some $(A,\vph)\in \A^\tp_{\mu,b}$.
On the other hand, by Proposition \ref{minuscule bijection} and \cite[Proposition 3.13]{Nie}, there exists a unique $\ld_\bl\in \cA_{\mu_\bl,b_\bl}^{\tp}$ with $\ld_1(1)+\cdots+\ld_1(n)=0$ such that $C=\pr(\overline{X_{\mu_\bl}^{\ld_\bl}(b_\bl)})$. 
Here $\pr\colon \cG r^d\rightarrow \cG r$ denotes the projection to the first factor.
The following theorem is established in \cite[Theorem 4.4]{Shimada3} by the author.
\begin{theo}
\label{constructionthm}
We have $\up_{\xi_j(\bb,\up)}=\up w_1^{-1}\cdots w_{j-1}^{-1}$ and $\xi_{\bl}(\bb,\up)\in \cA_{\mu_\bl, b_\bl}^\tp$.
If $\up'$ is an element in $\Upsilon(\bb)$ different from $\up'$, then $\xi_{\bl}(\bb, \up)\sim \xi_{\bl}(\bb, \up')$.
Let $\xi_\bl^0(\bb)$ be the unique cocharacter in $[\xi_\bl(\bb, \up)]$ such that $\xi_1^0(\bb)(1)+\cdots+\xi_1^0(\bb)(n)=0$.
Then for any $(A,\vph)\in \A^\tp_{\mu,b}$, there exists a unique $\bb\in \B_\mu(\ld_b)$ such that $\overline{S_{A,\vph}}=\pr(\overline{X_{\mu_\bl}^{\xi_\bl^0(\bb)}(b_\bl)})$.
\end{theo}
\begin{proof}
This is \cite[Theorem 4.4]{Shimada3}.
\end{proof}

This correspondence between $\A^\tp_{\mu,b}$ and $\B_\mu(\ld_b)$ above is compatible with the natural bijection in Chen-Zhu conjecture constructed by Nie in \cite{Nie}.

\begin{coro}
Let $(A,\vph)\in \A^\tp_{\mu,b}$. 
Let $\bb\in \B_\mu(\ld_b)$ such that $\overline{S_{A,\vph}}=\pr(\overline{X_{\mu_\bl}^{\xi_\bl^0(\bb)}(b_\bl)})$.
Then $(A,\vph)$ is cyclic if and only if $$\sum_{1\le j\le d}w_1^{-1}\cdots w_{j-1}^{-1}\wt(\bb_j)\in W_0\mu.$$
\end{coro}
\begin{proof}
By Lemma \ref{semi-module str}, we have $A=A^{\xi_1^0(\bb)}$.
Recall that $(A,\vph)$ is cyclic if and only if $b\xi_1^0(\bb)-\xi_1^0(\bb)\in W_0\mu$.
Since $b\xi_1^0(\bb)-\xi_1^0(\bb)$ is a conjugate of $b\xi_1(\bb,\up)-\xi_1(\bb,\up)$, this is also equivalent to $\up^{-1}b\xi_1(\bb,\up)-\up^{-1}\xi_1(\bb,\up)\in W_0\mu$.
By Theorem \ref{constructionthm}, 
$$\up^{-1}b\xi_1(\bb,\up)-\up^{-1}\xi_1(\bb,\up)=\sum_{1\le j\le d}w_1^{-1}\cdots w_{j-1}^{-1}\wt(\bb_j).$$
This finishes the proof.
\end{proof}

We call an element $\bb\in \B_\mu(\ld_b)$ cyclic if $$\ld(\bb)\coloneqq\sum_{1\le j\le d}w_1^{-1}\cdots w_{j-1}^{-1}\wt(\bb_j)\in W_0\mu.$$
Now we give another interpretation of Lemma \ref{dual cyclic}.
By Exmaple \ref{highest}, $\B_\mu^*$ is isomorphic to $\B_{\mu^*}$.
We denote by $\bb^*\in \B_{\mu^*}$ the dual of $\bb\in \B_\mu$.
Note that we have $(w\bb)^*=w\bb^*$ for any $w\in W_0$.
So if $\bb\in \B_\mu(\ld_b)$, then ${\bb^{\op}}^*=w_{\max}\bb^*\in\B_{\mu^*}(\ld_{b^*})$.
\begin{lemm}
\label{dual cyclic2}
We have $\ld({\bb^{\op}}^*)=-w(\bb)^{-1}\ld(\bb)+(d,\ldots, d)$.
In particular, $\bb\in \B_\mu(\ld_b)$ is cyclic if and only if ${\bb^{\op}}^*\in\B_{\mu^*}(\ld_{b^*})$ is cyclic.
\end{lemm}
\begin{proof}
Note that if $(\mu_1,\ldots,\mu_d)$ is Far-Eastern, then $(\mu_d^*,\ldots, \mu_1^*)$ is Far-Eastern.
So if we write $$\FE(\bb)=\bb_1\otimes\cdots \otimes \bb_d\quad \text{and}\quad \FE(\bb^{\op})=w_1\bb_1\otimes\cdots \otimes w_d\bb_d$$
in $\B_{\mu_1}\otimes\cdots\otimes\B_{\mu_d}$,
then we have $$\FE(\bb^*)=\bb_d^*\otimes\cdots \otimes \bb_1^*\quad \text{and}\quad \FE({\bb^{\op}}^*)=w_d\bb_d^*\otimes\cdots \otimes w_1\bb_1^*$$
in $\B_{\mu_d^*}\otimes\cdots\otimes \B_{\mu_1^*}$.
Thus $w({\bb^{\op}}^*)=w_d\cdots w_1=w(\bb)^{-1}, \Upsilon({\bb^{\op}}^*)=\Upsilon(\bb)$ and
\begin{align*}
\ld({\bb^{\op}}^*)&=\wt(w_d\bb_d^*)+w_d\wt(w_{d-1}\bb_{d-1}^*)+\cdots+w_d\cdots w_2\wt(w_1\bb_1^*)\\
&=-w(\bb)^{-1}\ld(\bb)+(d,\ldots, d),
\end{align*}
as desired.
\end{proof}

\subsection{Non-Cyclic Semi-standard Tableaux}
\label{NCtableaux}
The goal of this section is to specify the dominant cocharacters $\mu$ such that every $\bb\in \B_\mu(\ld_b)$ is cyclic.
Set $d=\mu(1)$.
\begin{lemm}
\label{dual inequality}
Assume that $n\geq 3$.
We have $d\geq 2\mn+\lfloor\frac{2m_0}{n}\rfloor+1$ or $d\geq 2\lfloor\frac{nd-m}{n}\rfloor+\lfloor\frac{2(n-m_0)}{n}\rfloor+1$.
\end{lemm}
\begin{proof}
It suffices to show that $d\le 2\mn+\lfloor\frac{2m_0}{n}\rfloor$ is equivalent to $d\geq 2\lfloor\frac{nd-m}{n}\rfloor+\lfloor\frac{2(n-m_0)}{n}\rfloor+1$.
Note that $\mn=\frac{m-m_0}{n}, \lfloor\frac{nd-m}{n}\rfloor=\frac{nd-m-(n-m_0)}{n}$.
So $d\le 2\mn+\lfloor\frac{2m_0}{n}\rfloor$ is equivalent to $(n-2)d\le 2(m-d-m_0)+n\lfloor\frac{2m_0}{n}\rfloor$, and $d\geq 2\lfloor\frac{nd-m}{n}\rfloor+\lfloor\frac{2(n-m_0)}{n}\rfloor+1$ is equivalent to $(n-2)d\le 2(m-d-m_0)+n(1-\lfloor\frac{2(n-m_0)}{n}\rfloor)$.
Then the assertion follows from the fact that $\lfloor\frac{2m_0}{n}\rfloor=0$ (resp.\ $1$) if and only if $\lfloor\frac{2(n-m_0)}{n}\rfloor=1$ (resp.\ $0$).
\end{proof}

\begin{lemm}
\label{mn=2}
Assume that $n\geq 3$.
Let $\mu\in \Y_+$ such that $d\geq 2\mn+\lfloor\frac{2m_0}{n}\rfloor+1, \mu(2)\geq 2$ and $\mn\geq 2$.
Then $\B_\mu(\ld_b)$ contains at least one non-cyclic element.
\end{lemm}
\begin{proof}
First we consider the case $n=3$.
In this case, we have $2\le \mu(2)\le \mn$.
Moreover, it is easy to check that $w(\bb)=s_2s_1$ and $s_1\in\supp(w_{d-\mn})$ for any $\bb\in \B_\mu(\ld_b)$.
Let $\bb$ be the unique element in $\B_\mu(\ld_b)$ whose second row contains exactly one \fbox{$3$}.
$$\begin{ytableau}
   1 & \cdots & \cdots & 1 & 2 & 3 & \cdots & 3 \\
   2 & \cdots & 2 & 3 
  \end{ytableau}
$$
Since $2\le \mu(2)\le \mn$, we have $$w_1^{-1}\cdots w_{d-\mu(2)}^{-1}\wt(\bb_{d-\mu(2)+1})=(0,1,1)\quad \text{and}\quad w_1^{-1}\cdots w_{d-1}^{-1}\wt(\bb_d)=(1,0,1).$$
Thus $\ld(\bb)\notin W_0\mu$ because $\mu(n)=0$.
This proves the case $n=3$.

In the rest of the proof, we assume that $n\geq 4$.
Let $\ld$ be a conjugate of $\ld_b$ such that $(\ld(1),\ld(2),\ld(3))=(\mn, \mn+\lfloor\frac{2m_0}{n}\rfloor, \mn+1)$ and $\ld(4)\geq \cdots \geq\ld(n)$.
Set $$\mu_0=(3\mn+\lfloor\frac{2m_0}{n}\rfloor+1-\min\{\mu(2),\mn\},\min\{\mu(2),\mn\},0,\ldots,0)\in \Y_+$$ and $\ld_0=(\ld(1),\ld(2),\ld(3),0,\ldots, 0)\in \Y$.
Note that we have $\mu(1)+\mu(2)\geq 3\mn+\lfloor\frac{2m_0}{n}\rfloor+1$.
Indeed if $\mu(1)+\mu(2)\le 3\mn+\lfloor\frac{2m_0}{n}\rfloor$, then by $\mu(1)\geq 2\mn+\lfloor\frac{2m_0}{n}\rfloor+1$, we have $\mu(2)\le \mn-1$.
This implies $\mu(3)+\cdots+\mu(n-1)\le (n-3)(\mn-1)$, or equivalently $3\mn+n+m_0-3\le \mu(1)+\mu(2)$, which is a contradiction.
Thus $Y_\mu$ contains $Y_{\mu_0}$.

Let $\bb_0$ be the unique element in $\B_{\mu_0}(\ld_0)$ whose second row contains exactly one \fbox{$3$}.
We will show that there exists $\bb'\in\B_\mu(\ld)$ that contains $\bb_0$.
It is easy to check that $\mu(n-1)\le \mn$ and $\mu(n-2)\le \mu_0(1)$.
So each column in $Y_{\mu/\mu_0}$ has at most $n-3$ boxes.
By filling each column with the numbers $1,\ldots, n-3$ so that the entries are starting with $1$ and increasing by one from top to bottom, we obtain a skew Young tableau of shape $Y_{\mu/\mu_0}$.
Let $k_i$ be the number of \fbox{$i$} in this tableau. Clearly we have $k_1\geq \cdots \geq k_{n-3}$.
$$\begin{ytableau}
   1 & \cdots & \cdots & 1 & 2 & 3 & \cdots & 3 &\null  \\
   2 & \cdots & 2 & 3 &\null &\null & \null \\
   \null &\null &\null &\null &\null  \\
   \null 
  \end{ytableau}
  \hspace{1cm}
  \begin{ytableau}
   \none & \none & \none & \none & \none & \none & \none & \none & 1 \\
   \none & \none & \none & \none & 1 & \cdots & 1 \\
   1 & \cdots & \cdots & 1 & 2  \\
   2 
  \end{ytableau}
$$
By $(\ld(4),\ldots,\ld(n))\pc (k_1,\ldots, k_{n-3})$ and Proposition \ref{Kostka}, there exists at least one skew Young tableau of shape $Y_{\mu/\mu_0}$ such that the number of \fbox{$i$} is $\ld(i+3)$ for each $1\le i\le n-3$.
By replacing $1,\ldots, n-3$ by $4,\ldots,n$ respectively, we obtain a skew Young tableau of shape $Y_{\mu/\mu_0}$ such that the number of \fbox{$i$} is $\ld(i)$ for each $4\le i\le n$.
Let $\bb'$ be the tableau obtained by joining $\bb_0$ and this skew tableau.
Clearly we have $\bb'\in \B_\mu(\ld)$, which shows our claim.

Let $\bb'\in \B_\mu(\ld)$ containing $\bb_0$, and let $\bb\in \B_\mu(\ld_b)$ be the conjugate of $\bb'$.
Then $s_2s_1\le w(\bb)$ and $s_1\in\supp(w_{d-\mn})$.
Let $k(\bb')$ be the number of \fbox{$4$} in the second row of $\bb'$.
If $k(\bb')<\mn$, then we have $$(w_1^{-1}\cdots w_{d-\min\{\mu(2),\mn\}}^{-1}\wt(\bb_{d-\min\{\mu(2),\mn\}+1}))(2)=1$$ and $$(w_1^{-1}\cdots w_{d-1}^{-1}\wt(\bb_d))(2)=0.$$
Thus $\ld(\bb)\notin W_0\mu$ and hence $\bb$ is non-cyclic.
If $k(\bb')\neq 0$, then $\ld(\bb)(1)=\mn-1$.

Assume that $\mu(3)<\mn-1$.
Then $\bb$ is always non-cyclic by the above argument.

Assume that $\mu(3)\geq 2$.
Let $\bb'_1=$\fbox{$j$} be the leftmost box in the third row of $\bb'$, and let $\bb'_2=$\fbox{$j'$} be the box right to $\bb'_1$.
Clearly $4\le j\le j'$.
$$
\begin{ytableau}
   1 & \cdots & \cdots & 1 & 2 & 3 & \cdots & 3 &\null  \\
   2 & \cdots & 2 & 3 & 4 & \cdots & 4 \\
   j & j' &\null &\null &\null  \\
   \null 
  \end{ytableau}
$$
Then in $\bb'$, all \fbox{$j-1$} are in the first or second row.
Since the number of \fbox{$j$} in the first or second row is less than $\wt(\bb')(j-1)$, there exists at least one \fbox{$j-1$} such that there is no box beneath it or the number in the box beneath it is greater than $j$.
So the tableau obtained by replacing $\bb'_1$ by the rightmost one among such \fbox{$j-1$} is semi-standard.
Repeating the same argument, we may assume $j=4$.
Similarly, if $\mn\geq 3$, we may also assume $j'=4$.
Indeed if $j'\geq 6$ and the leftmost column in $\bb'$ contains \fbox{$j'-1$} but does not contain \fbox{$j'$}, we replace $\bb'_2$ by this \fbox{$j'-1$}.
In other cases, by $\mn\geq 3$, there exists at least one \fbox{$j'-1$} such that there is no box beneath it or the number in the box beneath it is greater than $j'$, and we replace $\bb'_2$ by the rightmost \fbox{$j'-1$} among such \fbox{$j'-1$}.
Then the obtained tableau is semi-standard.
Thus if $\mn\geq 3$, there exists $\bb'$ containing $\bb_0$ such that $k(\bb')<\mn$, which is non-cyclic by the above argument.
If $\mn=2$ and $n=4$, then $\bb$ is non-cyclic because $k(\bb')<2$.
If $\mn=2$ and $n\geq 5$, we may also assume $j'=4$ and hence $\bb$ is non-cyclic unless the third row of $\bb'$ contains three \fbox{$5$}.
If $\mn=2, n\geq 5$ and the third row of $\bb'$ contains three \fbox{$5$}, then
$$(w_1^{-1}\cdots w_{d-2}^{-1}\wt(\bb_{d-1}))(4)=1\quad \text{and}\quad  (w_1^{-1}\cdots w_{d-1}^{-1}\wt(\bb_d))(4)=0.$$
Thus $\ld(\bb)\notin W_0\mu$ and hence $\bb$ is non-cyclic.
$$
\begin{ytableau}
   1 & 1 & 2 & 3 & 3 & \null & \null & \null &\null  \\
   2 & 3 & 4 & 4 & \null & \null & \null \\
   4 & 5 & 5 & 5 &\null  \\
   \null 
  \end{ytableau}
$$

Assume that $\mn=2$ and $\mu(3)=1$.
By the same argument as above, we may assume that the leftmost column of $\bb'$ contains \fbox{$4$}.
So $\bb$ is non-cyclic when $\ld(4)=2$.
If $\mu(1)> 5+\lfloor \frac{2m_0}{n}\rfloor$, we may assume that the first row of $\bb'$ also contains \fbox{$4$}.
This can be checked easily as above using $\mu(3)=1$.
Thus if $\mu(1)> 5+\lfloor \frac{2m_0}{n}\rfloor$, we obtain a non-cyclic $\bb$.
$$
\begin{ytableau}
   1 & 1 & 2  & 3 & 3 & 4 & \null & \null & \null   \\
   2 & 3 & 4  & \null &\null &\null & \null \\
   4 \\
   \null
  \end{ytableau}
$$
If $\mu(1)= 5+\lfloor \frac{2m_0}{n}\rfloor$, then we have $n=4$ or $5$. 
More precisely, we have $$\mu=(6,4,1,0), (5,5,1,1,0),(6,5,1,1,0),(6,6,1,0,0),\ \text{or}\ (6,6,1,1,0),$$
and $\bb'$ contains one of the following smaller Young tableaux when $\ld(4)=3$.
$$
\begin{ytableau}
   1 & 1 & 2 & 3 & 3  \\
   2 & 3 & 4 & 4  \\
   4 
  \end{ytableau}
\hspace{1cm}
\begin{ytableau}
   1 & 1 & 2 & 2 & 3 & 3  \\
   2 & 3 & 4 & 4  \\
   4 
  \end{ytableau}
$$
We can easily check that $\bb$ is non-cyclic in every case.

Putting things together, we have proved the lemma.
\end{proof}

\begin{lemm}
\label{mn=1}
Assume that $n\geq 4$.
Let $\mu\in \Y_+$ such that $d\geq 3+\lfloor\frac{2m_0}{n}\rfloor, \mu(2)\geq 2$ and $\mn=1$.
Then $\B_\mu(\ld_b)$ contains at least one non-cyclic element.
\end{lemm}
\begin{proof}
Let $\ld$ be a conjugate of $\ld_b$ such that $(\ld(1),\ld(2),\ld(3))=(\ld_b(1),\ld_b(2),\ld_b(3))$ and $\ld(4)\geq \cdots \geq\ld(n)$.
Assume that $(\ld_b(1),\ld_b(2),\ld_b(3))=(1,2,2)$ and $\mu(2)\geq 3$.
Similarly as the proof of Lemma \ref{mn=2}, we can easily show that there exists $\bb'\in \B_\mu(\ld)$ containing the following smaller Young tableau.
$$
  \begin{ytableau}
   1 & 2 & 3 & 4\\
   2 & 3 & 4
  \end{ytableau}
$$
Let $\bb\in \B_\mu(\ld_b)$ be the conjugate of $\bb'$.
If $\mu(3)<2$, then $\bb$ is non-cyclic because $\ld(\bb)(2)=2$.
If $\mu(3)\geq 2$, then similarly as the proof of Lemma \ref{mn=2}, we may assume that the second row of $\bb'$ does not contain \fbox{$5$}.
In this case, the conjugate $\bb\in \B_\mu(\ld_b)$ of $\bb'$ is non-cyclic because $$(w_1^{-1}\cdots w_{d-3}^{-1}\wt(\bb_{d-2}))(3)=1\quad\text{and}\quad  (w_1^{-1}\cdots w_{d-1}^{-1}\wt(\bb_d))(3)=0.$$

Assume that $(\ld_b(1),\ld_b(2),\ld_b(3))=(1,2,2)$ and $\mu(2)=2$.
Then there exists $\bb'\in \B_\mu(\ld)$ containing one of the following smaller Young tableaux.
$$
\begin{ytableau}
   1 & 2 & 3 &3 &4 \\
   2 & 4
  \end{ytableau}
  \hspace{1cm}
  \begin{ytableau}
   1 & 2 & 3 & 4\\
   2 & 4 \\
   3
  \end{ytableau}
$$
It is easy to check that the conjugate $\bb\in \B_\mu(\ld_b)$ of $\bb'$ is non-cyclic.

Assume that $(\ld_b(1),\ld_b(2),\ld_b(3))\neq(1,2,2)$.
Then there exists $\bb'\in \B_\mu(\ld)$ containing one of the following smaller Young tableaux.
$$
  \begin{ytableau}
   1 & 2 & 3 & 4 \\
   2 & 4 
  \end{ytableau}
  \hspace{1cm}
  \begin{ytableau}
   1 & 3 & 3 & 4 \\
   2 & 4 
  \end{ytableau}
  \hspace{1cm}
\begin{ytableau}
   1 & 3 & 4 \\
   2 & 4 
  \end{ytableau}
  \hspace{1cm}
$$
Let $\bb\in \B_\mu(\ld_b)$ be the conjugate of $\bb'$.
Since $\ld(\bb)(1)=1$, $\bb$ is non-cyclic if $\mu(3)=0$.
If $\mu(3)\geq 2$, then similarly as the proof of Lemma \ref{mn=2}, we may assume that the second row of $\bb'$ does not contain \fbox{$5$}.
In this case, $\bb$ is non-cyclic because 
$$(w_1^{-1}\cdots w_{d-2}^{-1}\wt(\bb_{d-1}))(3)=1\quad\text{and}\quad  (w_1^{-1}\cdots w_{d-1}^{-1}\wt(\bb_d))(3)=0.$$
If $\mu(3)=1$ and $\mu(1)>3+\lfloor\frac{2m_0}{n}\rfloor$, then we may also assume that the second row of $\bb'$ does not contain \fbox{$5$} and hence $\bb$ is non-cyclic.
If $\mu(3)=1$ and $\mu(1)=3+\lfloor\frac{2m_0}{n}\rfloor$, then we may assume that the leftmost column of $\bb'$ contains \fbox{$5$}.
We can easily check that $\bb$ is non-cyclic by an easy calculation.
$$
  \begin{ytableau}
   1 & 2 & 3 & 4   \\
   2 & 4 & 5 & \null  \\
   5 \\
   \null
  \end{ytableau}
  \hspace{1cm}
    \begin{ytableau}
   1 & 3 & 3 & 4   \\
   2 & 4 & 5 & \null  \\
   5 \\
   \null
  \end{ytableau}
  \hspace{1cm}
  \begin{ytableau}
   1 & 3 & 4    \\
   2 & 4 & 5  \\
   5 \\
   \null
  \end{ytableau}
$$
This finishes the proof.
\end{proof}

\begin{lemm}
\label{mn=0}
Assume that $n\geq 5$.
Let $\mu\in \Y_+$ such that $\mn=0$.
If $(1)\ \mu(2)\geq 2$ or $(2)\ d\geq 3, \mu(2)=1$,
then $\B_\mu(\ld_b)$ contains at least one non-cyclic element.
\end{lemm}
\begin{proof}
Let $1<i_1<i_2<\cdots<i_{m_0}=n$ be the integers such that $\ld_b(i_1)=\ld_b(i_2)=\cdots =\ld_b(i_{m_0})=1$.
Let $\bb$ be the Young tableau in $\B_\mu(\ld_b)$ obtained by filling $Y_\mu$ with $i_1,\ldots, i_{m_0}$ from top to bottom, starting from the leftmost column.
$$
\ytableausetup{boxsize=1.8em}
\begin{ytableau}
   i_1 & i_{k+1} & \cdots & i_{m}  \\
   i_2 & i_{k+2} & \vdots \\
   \vdots & \vdots \\
   i_{k}
  \end{ytableau}
$$
If (1) holds, then $\bb$ is non-cyclic because 
$$\wt(\bb_1)(i_m)=1\quad \text{and}\quad (w_1^{-1}\cdots w_{d-1}^{-1}\wt(\bb_d))(i_m)=0.$$
Let $k=\max\{i\mid\mu(i)\neq 0\}$.
If (2) holds, then the Young tableau $\bold c\in \B_\mu(\ld_b)$ obtained by replacing \fbox{$i_k$} by \fbox{$i_{k+1}$} in $\bb$ is non-cyclic because $\ld(\bold c)(i_k)=2$.
$$
\ytableausetup{boxsize=1.8em}
\begin{ytableau}
   i_1 & i_k & i_{k+2} & \cdots & i_{m}  \\
   \vdots   \\
   i_{k-1} \\
   i_{k+1}
  \end{ytableau}
$$
This finishes the proof.
\end{proof}

\begin{theo}
\label{cyclic classification}
Every $\bb\in \B_\mu(\ld_b)$ is cyclic if and only if $\mu$ has one of the following forms:
\begin{enumerate}[(i)]
\item $\omega_i$ with $1\le i\le n-1$ such that $i$ is coprime to $n$.
\item $\omega_1+\omega_i$ or $\omega_{n-1}+\omega_{n-i}$  with $1\le i\le n-1$ such that $i+1$ is coprime to $n$.
\item $(nr+i)\omega_1$ or $(nr+i)\omega_{n-1}$ with $r\geq 0$ and $1\le i\le n-1$ such that $i$ is coprime to $n$.
\item $(nr+i-j)\omega_1+\omega_j$ or $(nr+i-j)\omega_{n-1}+\omega_{n-j}$ with $r\geq 1$, $2\le j\le n-1$ and $1\le i\le n-1$ such that $i$ is coprime to $n$.
\end{enumerate}
\end{theo}
\begin{proof}
It is easy to check that every $\bb\in \B_\mu(\ld_b)$ is cyclic if $\mu$ is one of the cocharacters in (i), (ii), (iii) and (iv).
It remains to show that if $\mu$ does not belong to the list above, then $\B_\mu(\ld_b)$ contains at least one non-cyclic element.
By Lemma \ref{dual cyclic2} and Lemma \ref{dual inequality}, we may assume that $d\geq 2\mn+\lfloor\frac{2m_0}{n}\rfloor+1$.
Then this follows from Lemma \ref{mn=2}, Lemma \ref{mn=1} and Lemma \ref{mn=0}.
\end{proof}

\begin{rema}
Even if every top extended semi-module for $\mu$ is cyclic, there might be a non-cyclic extended semi-module for $\mu$.
In fact, such cases exist, see \S\ref{semi-module 12}.
\end{rema}

\section{The Semi-Module Stratification}
\label{the semi-module stratification}
Keep the notations and assumptions in \S3.

\subsection{The Semi-Module Stratification for $\omega_i$}
\label{semi-module i}
Recall that if $\mu$ is minuscule, then every extended semi-module is cyclic.
\begin{lemm}
\label{omega2}
For any $1\le j\le \frac{n-3}{2}(=\dim X_{\omega_2}(\tau^2))$, we have 
$$\A_{\omega_2,\tau^2}^j=
\begin{cases}
\{[\chi_{2,n-1}^\vee+\chi_{4,n-3}^\vee+\cdots+\chi_{j,n-j+1}^\vee]\} & (j\ \text{even})\\
\{[\chi_{1,n}^\vee+\chi_{3,n-2}^\vee+\cdots+\chi_{j,n-j+1}^\vee]\} & (j\ \text{odd}).
\end{cases}$$
\end{lemm}
\begin{proof}
By (the proof of) \cite[Proposition 5.5]{Viehmann2}, each normalized semi-module for $2,n$ is of the form $A_j=(2\N-j)\cup (\N+j+1)$ for some $1\le j\le \frac{n-3}{2}$.
It is easy to check that
$$A_j=
\begin{cases}
A^{\chi_{2,n-1}^\vee+\chi_{4,n-3}^\vee+\cdots+\chi_{j,n-j+1}^\vee} & (j\ \text{even})\\
A^{\chi_{1,n}^\vee+\chi_{3,n-2}^\vee+\cdots+\chi_{j,n-j+1}^\vee} & (j\ \text{odd}).
\end{cases}$$
Let $(A_j,\vph_j)$ be the cyclic semi-module for $\omega_2$.
Then $n-2-j,n-1+j\in \bar{A_j}$ and $\vph_j(n-2-j)=\vph_j(n-1+j)=1$.
It is also easy to check that $|\cV(A_j,\vph_j)|=j$.
This finishes the proof.
\end{proof}

\begin{lemm}
\label{omega3}
Assume that $n=7$.
Then $\dim X_{\omega_3}(\tau^3)=3$ and
\begin{align*}
\A_{\omega_3,\tau^3}^1=\{[\chi_{1,7}^\vee]\},\quad
\A_{\omega_3,\tau^3}^2=\{[\chi_{1,6}^\vee],[\chi_{2,7}^\vee]\},\quad
\A_{\omega_3,\tau^3}^3=\{[\chi_{3,5}^\vee]\}.
\end{align*}

Assume that $n=8$.
Then $\dim X_{\omega_3}(\tau^3)=4$ and
\begin{align*}
\A_{\omega_3,\tau^3}^1=\{[\chi_{1,8}^\vee]\},\quad &\A_{\omega_3,\tau^3}^2=\{[\chi_{1,7}^\vee],[\chi_{2,8}^\vee]\},\\
\A_{\omega_3,\tau^3}^3=\{[\chi_{2,6}^\vee],[\chi_{3,7}^\vee]\},\quad &\A_{\omega_3,\tau^3}^4=\{[\chi_{1,8}^\vee+\chi_{4,5}^\vee]\}.
\end{align*}
\end{lemm}
\begin{proof}
Using Lemma \ref{type}, we can easily check the lemma by an easy calculation.
\end{proof}

\subsection{The Semi-Module Stratification for $\omega_1+\omega_{n-2}$}
\label{semi-module 1n-2}
Throughout this subsection, we set $\mu=\omega_1+\omega_{n-2}$.
Also we assume that $n\geq 4$.
\begin{lemm}
\label{semi-module 1n-2 lemm}
Every extended semi-module for $\mu$ is cyclic.
For any $0\le j\le n-2(=\dim X_{\mu}(b))$, we define $\A_{\mu,b}^j$ similarly as in \S \ref{minuscule case}.
Then we have $\A_{\mu,b}^0=\emptyset$ and $|\A_{\mu,b}^j|=j$.
More precisely, if $j$ is odd, then $\A_{\mu,b}^j$ is equal to
\begin{align*}
\{[\chi_{1,n-j+1}^\vee],&[\chi_{1,n-j+3}^\vee+\chi_{2,n-j+2}^\vee],\ldots,\\
&[\chi_{1,n}^\vee+\chi_{2,n-1}^\vee+\cdots+\chi_{\frac{j+1}{2},n-\frac{j-1}{2}}^\vee],\ldots,
[\chi_{j-2,n}^\vee+\chi_{j-1,n-1}^\vee],[\chi_{j,n}^\vee]\},
\end{align*}
and if $j$ is even, then $\A_{\mu,b}^j$ is equal to
\begin{align*}
\{[\chi_{1,n-j+1}^\vee],&[\chi_{1,n-j+3}^\vee+\chi_{2,n-j+2}^\vee],\ldots,\\
&[\chi_{1,n-1}^\vee+\chi_{2,n-2}^\vee+\cdots+\chi_{\frac{j}{2},n-\frac{j}{2}}^\vee],\ldots,
[\chi_{j-2,n}^\vee+\chi_{j-1,n-1}^\vee],[\chi_{j,n}^\vee]\}.
\end{align*}
\end{lemm}
\begin{proof}
Let $(A,\vph)$ be an extended semi-module for $\mu$.
Let $\mu'$ be the type of $A$.
If $(A,\vph)$ is non-cyclic, then by Lemma \ref{cyclic}, $\mu'_{\dom}\prec \mu$, i.e., $\mu'_{\dom}=\omega_{n-1}$. 
By Lemma \ref{type}, we have $A=\{0,1,\ldots,n-1,\ldots\}$.
By Definition \ref{exsemidef} (3), $\vph(a)=\max\{k\mid a+n-1-kn\in A\}$ for all $a\in A$.
This contradicts to the assumption that $(A,\vph)$ is non-cyclic.
Thus $(A,\vph)$ is cyclic.

Since $\mu'$ satisfies $\nu_b\pc w_{\max}\mu'$, it is easy to check that
$$w_{\max}\mu'=s_{l+1}\cdots s_{n-3}s_{n-2}s_{k-1}\cdots s_2s_1\mu$$
for some $1\le k\le n-2$ and $k\le l\le n-2$.
Let $\bar A=\{a_0,a_1\ldots, a_{n-1}\}$ with $a_0=\min \bar A$.
Then we have $\vph(a_0)=0,\vph(a_{n-l-1})=0,\vph(a_{n-k})=2$ and $\vph(a_i)=1$ for $i\neq 0,n-l-1,n-k$.
Thus 
\begin{align*}
\cV(A,\vph)=&\{(a_{n-k},a_{n-l-1}+n),(a_{n-k},a_{n-l}),(a_{n-k},a_{n-l+1}),\ldots,(a_{n-k},a_{n-k-1})\}\\
\cup&\{(a_{n-k+1},a_{n-l-1}),(a_{n-k+2},a_{n-l-1}),\ldots, (a_{n-1},a_{n-l-1})\}
\end{align*}
and $|\cV(A,\vph)|=l$.
Then by Proposition \ref{UAphi}, the description of $\A_{\mu,b}^l$ for each $l$ in the lemma follows from direct computation.
\end{proof}

\subsection{The Semi-Module Stratification for $\omega_1+\omega_{n-3}$}
\label{semi-module 1n-3}
Throughout this subsection, we set $\mu=\omega_1+\omega_{n-3}$.
Also we assume that $n\geq 7$.
\begin{lemm}
\label{semi-module 1n-3 lemm}
Every extended semi-module for $\mu$ is cyclic.
For any $1\le j\le \frac{3n-9}{2}(=\dim X_{\mu}(b))$,  we define $\A_{\mu,b}^j$ similarly as in \S \ref{minuscule case}.
Then $|\A_{\mu,b}^{\frac{3n-9}{2}}|=n-3$ and $|\A_{\mu,b}^{\frac{3n-11}{2}}|\le 2(n-4)$.
\end{lemm}
\begin{proof}
Using Lemma \ref{omega2}, we can show the first assertion similarly as the proof of Lemma \ref{semi-module 1n-2 lemm}.
Indeed, for any semi-module $A^\ld$ in Lemma \ref{omega2}, there exists a unique $\vph$ such that $(A^\ld,\vph)$ is an extended semi-module for some $\mu\in \Y_+$.
The equality $|\A_{\mu, b}^{\frac{3n-9}{2}}|=n-3$ follows from the Chen-Zhu conjecture.

Let $(A,\vph)$ be an extended semi-module for $\mu$ with type $\mu'(\in W_0\mu)$.
Let $0<k_1<k_2$ be integers such that $\mu'(1)=\mu'(k_1+1)=\mu'(k_2+1)=0$, and let $l$ be an integer such that $\mu'(l+1)=2$.
Assume that $\nu_b\pc w_{\max}s_{k_2+1}\mu'$.
Let $(B,\psi)$ be an extended semi-module for $\mu$ with type $s_{k_2+1}\mu'$.
Let $a_0=\min \bar A$ (resp.\ $b_0=\min \bar B$) and let inductively $a_i=a_{i-1}+n-2-\mu'(i)n$ (resp.\ $b_i=b_{i-1}+n-2-(s_{k_2+1}\mu')(i)n$) for $i=1,\ldots, n$.
Then $a_0=a_n$ (resp.\ $b_0=b_n$) and $\{a_0,a_1,\ldots,a_{n-1}\}=\bar A$ (resp.\  $\{b_0,b_1,\ldots,b_{n-1}\}=\bar B$).
We will show that if $l>k_2+1$ (resp.\ $l=k_2+1$), then $|\cV(B,\psi)|\le|\cV(A,\vph)|$ (resp.\ $|\cV(B,\psi)|<|\cV(A,\vph)|-1$).
Moreover, the equality does not hold if $k_2-k_1\le 3.$

Note that we have $\vph(a_{0})=\vph(a_{k_1})=\vph(a_{k_2})=0,\vph(a_l)=2,\psi(b_0)=\psi(b_{k_1})=\psi(b_{k_2+1})=0,\psi(b_l)=2$.
Note also that
\begin{align*}
\cV(A,\vph)=&\{(a,a')\mid \text{$a\in \bar A$ with $\vph(a)=1$, $a'=a_{k_1}$ or $a_{k_2}$}\}\\
\sqcup &\{(a_l,a')\mid a_l<a',\vph(a')<2\}
\end{align*}
and
\begin{align*}
\cV(B,\psi)=&\{(b,b')\mid \text{$b\in \bar B$ with $\psi(b)=1$, $b'=b_{k_1}$ or $b_{k_2+1}$}\}\\
\sqcup &\{(b_l,b')\mid b_l<b',\psi(b')<2\}.
\end{align*}
Let $\cV(A,\vph)_1$ (resp.\ $\cV(B,\psi)_1$) be the first subset in $\cV(A,\vph)$ (resp.\ $\cV(B,\psi)$) above, and let $\cV(A,\vph)_2$ (resp.\ $\cV(B,\psi)_2$) be its complement.

If $l>k_2+1$, then it follows that 
$$b_k=\begin{cases}
a_k+1 & (k\neq k_2+1)\\
a_k+1-n & (k=k_2+1)
\end{cases},\quad
\psi(b_k)=\begin{cases}
\vph(a_k) & (k\neq k_2,k_2+1)\\
1-\vph(a_k) & (k=k_2,k_2+1).
\end{cases}$$
In particular, $b_{k_2+1}-1=a_{k_2}-2$.
So $|\cV(B,\psi)_1|>|\cV(A,\vph)_1|$ implies that $|\cV(B,\psi)_1|=|\cV(A,\vph)_1|+1$ and $b_{k_2}<b_{k_1}$.
By the fact $(a_l, a_{k_2+1})\in \cV(A,\vph)_2$, we always have $|\cV(B,\psi)_2|< |\cV(A,\vph)_2|$.
Thus $|\cV(B,\psi)|\le|\cV(A,\vph)|$.
Moreover, if $k_2-k_1\le 3$, then the equality does not hold because $b_{k_2}\geq b_{k_1}$.

If $l=k_2+1$, then it follows that 
$$b_k=\begin{cases}
a_k+2 & (k\neq k_2+1)\\
a_k+2-2n & (k=k_2+1)
\end{cases},\quad
\psi(b_k)=\begin{cases}
\vph(a_k) & (k\neq k_2,k_2+1)\\
2-\vph(a_k) & (k=k_2,k_2+1).
\end{cases}$$
In particular, $b_{k_2+1}-2=a_{k_2}-2-n$.
By $\nu_b\pc w_{\max}s_{k_2+1}\mu'$, we have $k_2\le \frac{n-3}{2}$.
Using this, we can easily check that $|\cV(B,\psi)|<|\cV(A,\vph)_1|$ and $\cV(A,\vph)_2=\{(a_{k_2+1},a_{k_2}+n)\}$.
Thus $|\cV(B,\psi)|<|\cV(A,\vph)|-1$.

Assume that $\nu_b\pc w_{\max}s_{k_1+1}\mu'$.
Let $(C,\chi)$ be an extended semi-module for $\mu$ with type $s_{k_1+1}\mu'$.
Similarly as above, we can show that if $l\geq k_1+1$, then $|\cV(C,\chi)|\le |\cV(A,\vph)|$.
Therefore, $|\cV(A,\vph)|\geq\frac{3n-11}{2}$ holds only if $k_2=2$ or $l>k_2=3$.
From this and $|\A_{\mu,b}^{\frac{3n-9}{2}}|=n-3$, we obtain $|\A_{\mu,b}^{\frac{3n-11}{2}}|\le 2(n-4)$.
\end{proof}

\subsection{The Semi-Module Stratification for $\omega_1+\omega_2,\omega_{n-3}+\omega_{n-1}$}
\label{semi-module 12}

\begin{lemm}
\label{semi-module 12 n=5}
Assume that $n=5$.
Set $\mu=\omega_1+\omega_2$.
Then every extended semi-module for $\mu$ is cyclic.
For any $1\le j\le 3(=\dim X_{\mu}(b))$, we define $\A_{\mu,b}^j$ similarly as in \S \ref{minuscule case}.
Then
\begin{align*}
\A_{\mu,b}^0=\emptyset,
\A_{\mu,b}^1=\emptyset,
\A_{\mu,b}^2=\{\chi_{1,4}^\vee,\chi_{2,5}^\vee\},
\A_{\mu,b}^3=\{\chi_{2,3}^\vee,\chi_{3,4}^\vee\}.
\end{align*}
\end{lemm}
\begin{proof}
The first assertion follows similarly as the proof of Lemma \ref{semi-module 1n-2 lemm}.
The second assertion follows from direct computation.
\end{proof}

\begin{lemm}
\label{semi-module 12 n=7,8}
Assume that $n=7$ or $8$.
Let $\mu$ be $\omega_1+\omega_2$ or $\omega_4+\omega_{n-1}$.
Then there exists a non-cyclic extended semi-module for $\mu$.
\end{lemm}
\begin{proof}
As described in Lemma \ref{omega3}, there exists a unique top cyclic extended semi-module $(A^\ld,\vph)$ for $\omega_3$.
We define $\vph'\colon \Z\rightarrow \N\cup \{-\infty\}$ by setting
$$\vph'(a)=
\begin{cases}
\vph(a) & (a\neq 1)\\
0 & (a=1).
\end{cases}$$
Then it is straightforward to check that $(A^\ld,\vph')$ is a non-cyclic extended semi-module for $\omega_1+\omega_2$.
The proof for $\omega_4+\omega_{n-1}$ is similar.
\end{proof}

\section{The Ekedahl-Oort Stratification}
\label{the Ekedahl-Oort stratification}
Keep the notations and assumptions in \S3.
For $\mu\in \Y_+$, set $$\SAdm(\mu)_{\cyc}=\{w\in \SAdm(\mu)\mid \text{$p(w)$ is $n$-cycle}\}.$$
By Theorem \ref{empty}, $X_w(b)\neq \emptyset$ if $w\in \SAdm(\mu)_{\cyc}$.
\subsection{The Ekedahl-Oort Stratification for $\mu=\omega_i$}
\label{EO i}
Throughout this subsection, we set $\mu=\omega_i$ and $c=s_is_{i+1}\cdots s_{n-1}s_{i-1}\cdots s_2s_1$.
By \cite[Theorem 2.7]{HNY}, we have $\dim X_{\vp^\mu c}(b)=\dim X_\mu(b)=\la\mu,\rho\ra-\frac{n-1}{2}$.

Let $w\in \tW$.
There exists a positive integer $k$ such that $w^k=\vp^\ld$ for some $\ld\in \Y$.
Then we set $\nu_w=\ld/k \in X_*(T)_{\Q}$.
This is independent of the choice of $k$.
Clearly $\supp_\sigma(w)\neq \tS$ implies that $\nu_w$ is central (cf.\ \cite[Lemma 1.1]{Lim23}).
\begin{lemm}
\label{n=9}
Assume that $n\geq 9$ and $4\le i\le n-4$.
Set $y=cs_is_{i+1}s_{i-1}=(1\ i+1\ i+3\ i+4\ \cdots\ n\ i\ i-2\ \cdots\  3\ 2)(i-1\ i+2)$.
Then we have $\vp^\mu y\in \SAdm(\mu)$ and $X_{\vp^\mu y}(b)\neq \emptyset$.
\end{lemm}
\begin{proof}
Under the assumption in the lemma, we have $\ell(\vp^\mu y)=\la\mu, 2\rho\ra-\ell(y)$ and hence $\vp^\mu y\in \SAdm(\mu)$ (cf.\ \cite[(2.4.5)]{Macdonald}).
Note that $\supp_\sigma(\vp^\mu y)=\tS$ because $\nu_{\vp^\mu y}$ is not central.
So, by Lemma \ref{LPr} and Theorem \ref{empty}, $X_{\vp^\mu y}(b)\neq \emptyset$ is equivalent to saying $\supp(ryr^{-1})\subsetneq S$ for any $r\in W_0$ such that $r(\Phi_+\setminus \Phi_{\vp^\mu y})\subset \Phi_+$.
It is easy to check that
$$\Phi_{\vp^\mu y}=\Phi_{\{\chi_{1,2},\chi_{2,3},\ldots,\chi_{i,i+1}\}}\cup \Phi_{\{\chi_{i,i+1},\chi_{i+1,i+2},\ldots,\chi_{n-1,n}\}}\cup\{\chi_{i-2,i+2},\chi_{i-1,i+2},\chi_{i-1,i+3}\}.$$
In particular, we have $\chi_{1,i+2},\chi_{i-1,n}\in \Phi_+\setminus \Phi_{\vp^\mu y}$.
Note that we can decompose $ryr^{-1}$ into disjoint cycles as $$(r(1)\ r(i+1)\ r(i+3)\ r(i+4)\ \cdots\ r(n)\ r(i)\ r(i-2)\ \cdots\  r(3)\ r(2))(r(i-1)\ r(i+2))$$
for any $r\in W_0$.
So if $ryr^{-1}\in \bigcup_{J\subsetneq S}W_J$, then $(r(i-1)\ r(i+2))=(1\ 2)$ or $(n-1\ n)$.
This implies that $r\chi_{1,i+2}$ or $r\chi_{i-1,n}$ is negative and hence that $r$ does not satisfy $r(\Phi_+\setminus \Phi_{\vp^\mu y})\subset \Phi_+$.
Thus we have $X_{\vp^\mu y}(b)\neq \emptyset$.
\end{proof}

\begin{lemm}
\label{n=10}
Assume that $n\geq 10$ and $i=3$ (resp.\ $i=n-3$).
Set $y=cs_3s_4s_5s_6s_2$ (resp.\ $y=cs_{n-3}s_{n-4}s_{n-5}s_{n-6}s_{n-2}$).
Then we have $\vp^\mu y\in \SAdm(\mu)$ and $X_{\vp^\mu y}(b)\neq \emptyset$.
\end{lemm}
\begin{proof}
We only treat the case $i=3$.
The proof for the case $i=n-3$ is similar.

The first assertion is easy.
To show the second assertion,
by Lemma \ref{LPr} and Theorem \ref{empty}, it suffices to check that $ryr^{-1}\notin \bigcup_{J\subsetneq S}W_J$ for any $r\in W_0$ such that $r(\Phi_+\setminus \Phi_{\vp^\mu y})\subset \Phi_+$.
By an explicit calculation, it follows that $\chi_{1,7},\chi_{2,9}\in \Phi_+\setminus \Phi_{\vp^\mu y}$ and
\begin{align*}
ryr^{-1}=(r(1)\ r(4)\ r(6)\ r(8)\ r(9)\ \cdots\ r(n)\ r(3))(r(2)\ r(5)\ r(7)).
\end{align*}
If $ryr^{-1}\in \bigcup_{J\subsetneq S}W_J$, then $(r(2)\ r(5)\ r(7))$ is equal to $(1\ 2\ 3)$ or $(n-2\ n-1\ n)$.
This implies that $r$ does not satisfy $r(\Phi_+\setminus \Phi_{\vp^\mu y})\subset \Phi_+$.
Thus we have $X_{\vp^\mu y}(b)\neq \emptyset$.
\end{proof}

\begin{lemm}
\label{n=10 Coxeter}
Assume that $n\geq 10$ and $i=3$ (resp.\ $i=n-3$).
Let $y$ be $cs_is_{i-1}$ or $cs_is_{i+1}$.
Then we have $\vp^\mu y\in \SAdm(\mu)$ and $X_{\vp^\mu y}(b)\neq \emptyset$.
\end{lemm}
\begin{proof}
The proof is similar to the proof of Lemma \ref{n=9} and Lemma \ref{n=10}.
Note that $y$ is $n$-cycle in this case.
\end{proof}

\begin{prop}
\label{EO i not}
Assume that $n\geq 9$ and $3\le i\le n-3$.
Then the semi-module stratification of $X_\mu(b)$ is not a refinement of the Ekedahl-Oort stratification.
\end{prop}
\begin{proof}
First assume that $n\geq 9$ and $4\le i\le n-4$.
Let $\vp^\mu y\in \SW$ be as in Lemma \ref{n=9}.
Let $\cT$ be a reduction tree of $\vp^\mu y$.
By Proposition \ref{rational point}, we have 
$$|X_{\vp^\mu y}(b)^{0,\sigma}|=\sum_{\unp}(q-1)^{\ell_I(\unp)}q^{\ell_{II}(\unp)},$$
where $\unp$ runs over all the reduction paths in $\cT$ with $\en(\unp)=\tau^m$.
Set $d=\dim X_\mu(b)=\la\mu,\rho\ra-\frac{n-1}{2}$.
Suppose that the semi-module stratification of $X_\mu(b)$ is a refinement of the Ekedahl-Oort stratification.
By Lemma \ref{injective}, Proposition \ref{DL method prop} and $\dim X_{\vp^\mu c}(b)=d$, we have $\ell_I(\unp)+\ell_{II}(\unp)\le \dim X_{\vp^\mu y}(b)\le d-1$ for any $\unp$.
On the other hand, we have $\ell_I(\unp)+2\ell_{II}(\unp)=\ell(\vp^\mu y)=2d-3$.
Thus we have $\ell_I(\unp)+\ell_{II}(\unp)=d-1$ and $\ell_I(\unp)=1$ for any $\unp$.
It follows that
$$|X_{\vp^\mu y}(b)^{0,\sigma}|=k(q-1)q^{d-2},$$
where $k\geq 1$ is the number of irreducible components of $X_{\vp^\mu y}(b)^0$.
Again by Lemma \ref{injective} and the fact that each $S_{A,\vph}$ is locally closed, we have $|\{(A,\vph)\mid \dim S_{A,\vph}=d-1, S_{A,\vph}\subseteq \pi(X_{\vp^\mu y}(b)^0)\}|=k$.
By Lemma \ref{cyclic}, it follows that 
$$k(q-1)q^{d-2}=kq^{d-1}+\sum_{0\le j\le d-2}k_j q^j,\qquad k_j\geq 0,$$
which is a contradiction.
This implies the proposition in this case.

Next assume that $n\geq 10$ and $i=3, n-3$.
Let $\vp^\mu y\in \SW$ be as in Lemma \ref{n=10}.
Suppose that the semi-module stratification of $X_\mu(b)$ is a refinement of the Ekedahl-Oort stratification.
Similarly as above, we can check that
$$\dim X_{cs_is_{i-1}}(b)=X_{cs_is_{i+1}}(b)=d-1.$$
By Lemma \ref{injective} and Proposition \ref{semi-module number}, we have $\dim X_{\vp^\mu y}(b)\le d-2$.
Similarly as above, it follows that 
$$k(q-1)q^{d-3}=kq^{d-2}+\sum_{0\le j\le d-3}k_j q^j,\qquad k_j\geq 0,$$
where $k$ is the number of irreducible components of $X_{\vp^\mu y}(b)^0$.
This is a contradiction, which finishes the proof.
\end{proof}

The following proposition is the complement of Proposition \ref{EO i not}.
\begin{prop}
\label{EO i refinement}
We have
\begin{align*}
\SAdm(\omega_1)_{\cyc}&=\{\tau\},\\
\SAdm(\omega_2)_{\cyc}&=\{\tau^2, s_0s_{n-1}\tau^2,s_0s_{n-1}s_{n-2}s_{n-3}\tau^2,\ldots, s_0s_{n-1}\cdots s_5s_4\tau^2\} &(n\geq 5),\\
\SAdm(\omega_3)_{\cyc}&=\{\tau^3,s_0s_6\tau^3,s_0s_6s_1s_0\tau^3,s_0s_6s_5s_1\tau^3,s_0s_6s_5s_1s_0s_6\tau^3\} &(n=7),\\
\SAdm(\omega_3)_{\cyc}&=\{\tau^3,s_0s_1\tau^3,s_0s_7s_6s_5\tau^3,s_0s_7s_6s_1\tau^3,s_0s_7s_6s_5s_1s_0\tau^3,\\
&\hspace{3.94cm}s_0s_7s_6s_1s_0s_7\tau^3,s_0s_7s_6s_5s_1s_0s_7s_6\tau^3\} &(n=8).
\end{align*}
Let $\vp^\mu y\in \SW$ be one of the elements above.
Then there exists $v\in \LP(\vp^\mu y)$ such that $v^{-1}yv$ is a Coxeter element.
Moreover, $X_w(b)=\emptyset$ for any $w\in \SAdm(\mu)\setminus \SAdm(\mu)_{\cyc}$, and the semi-module stratification of $X_{\mu}(b)$ is a refinement of the Ekedahl-Oort stratification.
\end{prop}
\begin{proof}
The equalities in the proposition follow from easy calculations.
For other statements, we only prove the case for $\omega_2$.
Other cases can be checked similarly.

Set $d=\frac{n-3}{2}$.
For $0\le j\le d$, we set $w_j=s_0s_{n-1}\cdots s_{n-2j+1}\tau^2$.
Then $\ell(w_j)=2j$ and $$p(w_j)=(1\ 3\ 5\ \cdots\  n-2j\ n-2j+1\ \cdots\ n\ 2\ 4\ \cdots\ n-2j-1).$$
Also it is easy to check that $$\Phi_+\setminus \Phi_{w_j}=\{\chi_{1,n-2j+1},\ldots,\chi_{1,n-1},\chi_{1,n}\}.$$
Clearly there exists $r\in W_0$ with $r(\Phi_+\setminus \Phi_{w_j})\subset \Phi_+$ such that $rp(w_j)r^{-1}$ is a Coxeter element (cf.\ \cite[Lemma 5.1]{Shimada3}).

For an integer $j$, let $0\le [j]<n$ denote its residue modulo $n$.
For $a,b\in \N$ with $a-b\in 2\Z$, we define
$t_{a,b}=s_{[b-2]}\cdots s_{[a+2]}s_{[a]}$.
Set 
\begin{align*}
w_{j,0}=w_j,w_{j,1}=t_{0,n-2j+1}w_jt_{0,n-2j+1}^{-1}, w_{j,2}=t_{n-1,n-2j+2}t_{0,n-2j+1}w_jt_{0,n-2j+1}^{-1}t_{n-1,n-2j+2}^{-1},\\
\ldots, w_{j,j}=t_{n-j+1,n-j}\cdots t_{n-1,n-2j+2}t_{0,n-2j+1}w_jt_{0,n-2j+1}^{-1}t_{n-1,n-2j+2}^{-1}\cdots t_{n-j+1,n-j}^{-1}.
\end{align*}
It is easy to check that the simple reflections in $t_{0,n-2j+1},t_{n-1,n-2j+2},\ldots,t_{n-j+1,n-j}$ define
\begin{align*}
w_j=w_{j,0}\rightarrow_{\sigma}w_{j,1}=s_{n-1}s_{n-2}\cdots s_{n-2j+2}\tau^2&\rightarrow_{\sigma}w_{j,2}=s_{n-2}s_{n-3}\cdots s_{n-2j+3}\tau^2\\
&\rightarrow_{\sigma}\cdots \rightarrow_{\sigma} w_{j,j}=\tau^2.
\end{align*}
Let $\unp_j$ be the reduction path (in a suitable reduction tree) defined by this reduction.
By Lemma \ref{injective}, Proposition \ref{decomposition}, Proposition \ref{rational point} and Proposition \ref{semi-module number}, it follows that $X_{w_j}(\tau^2)=X_{\unp_j}$ and $X_w(\tau^2)=\emptyset$ for any $w\in \SAdm(\omega_2)\setminus \SAdm(\omega_2)_{\cyc}$.
Note that $X_{\tau^2}(\tau^2)^0=\{I\}$.
It is easy to check that $$\ell(t_{n-j+1,n-j}\cdots t_{n-1,n-2j+2}t_{0,n-2j+1})=\ell(t_{n-j+1,n-j})+\cdots+ \ell(t_{n-1,n-2j+2})+\ell(t_{0,n-2j+1}).$$
Thus by Proposition \ref{DL method prop} (cf.\ \cite[\S 3.3]{Shimada2}), each element $gI$ in $X_{w_j}(\tau^2)^0$ is contained in a Schubert cell associated to $t_{n-j+1,n-j}\cdots t_{n-1,n-2j+2}t_{0,n-2j+1}$.
By Lemma \ref{omega2}, it follows that $\pi(X_{w_j}(b)^0)$ is equal to the unique semi-module stratum of dimension $j$.
This shows that the semi-module stratification of $X_{\mu}(b)$ is a refinement of the Ekedahl-Oort stratification.
\end{proof}

\subsection{The Ekedahl-Oort Stratification for $\omega_1+\omega_{n-2}$}
\label{EO 1n-2}
Throughout this subsection, we set $\mu=\omega_1+\omega_{n-2}$.
Also we assume that $n\geq 4$.
Note that the unique dominant cocharacter $\mu'$ with $\mu'\prec \mu$ is $\omega_{n-1}$.
Clearly we have $\SAdm(\omega_{n-1})_{\cyc}=\{\tau^{n-1}\}$ and the semi-module stratification of $X_{\omega_{n-1}}(\tau^{n-1})$ is a refinement of the Ekedahl-Oort stratification.
\begin{prop}
\label{EO 1n-2 refinement}
For any $1\le j\le n-2(=\dim X_\mu(b))$, there exist exactly $j$ elements of length $2j$ in $\SAdm(\mu)_{\cyc}^\circ\coloneqq\SAdm(\mu)_{\cyc}\setminus \{\tau^{n-1}\}$.
Let $\vp^\mu y\in \SW$ be one of such elements.
Then there exists $v\in \LP(\vp^\mu y)$ such that $v^{-1}yv$ is a Coxeter element.
Moreover, $X_w(b)=\emptyset$ for any $w\in \SAdm(\mu)\setminus \SAdm(\mu)_{\cyc}$, and the semi-module stratification of $X_{\mu}(b)$ is a refinement of the Ekedahl-Oort stratification.
\end{prop}
\begin{proof}
We first prove by induction on $n$ that there exist at least $j$ elements of length $2j$ in $\SAdm(\mu)_{\cyc}^\circ$, each of which has finite part $y$ such that $ryr^{-1}$ is a Coxeter element for some $r\in W_{\{s_2,\ldots,s_{n-2}\}}$ satisfying $r(\Phi_+\setminus \Phi_{\vp^\mu y})\subset \Phi_+$ (cf.\ Lemma \ref{LPr}).
Note that if $y\in W_0$ satisfies 
\begin{align*}
\text{$y^{-1}(2)<y^{-1}(3)<\cdots <y^{-1}(n-2)$ and $y^{-1}(n-1)<y^{-1}(n)$,} \tag{$\ast$}
\end{align*}
then by \cite[Lemma 4.4]{Shimada2}, we have $\vp^\mu y\in \SAdm(\mu)$.
In particular, since $\ell(\vp^\mu)=3n-5$, $\vp^\mu y$ is an element of length $2j$ in $\SAdm(\mu)_{\cyc}^\circ$ for any $n$-cycle $y$ of length $3n-2j-5$.
If $n=4$, then $s_1s_2s_3,s_2s_3s_1$ and $s_1s_2s_3s_1s_2$ are $4$-cycles satisfying $(\ast)$.
Moreover, $s_2(s_1s_2s_3s_1s_2)s_2=s_1s_2s_3$ is a Coxeter element and $s_2(\Phi_+\setminus \Phi_{\vp^\mu s_1s_2s_3s_1s_2})\subset \Phi_+$.
So the claim is true for $n=4$.

Suppose that $n\geq 5$ and the claim is true for $n-1$.
Let $y$ be a $(n-1)$-cycle in $W_{\{s_1,s_2,\ldots,s_{n-2}\}}$ such that
$y^{-1}(2)<y^{-1}(3)<\cdots <y^{-1}(n-3)$ and $y^{-1}(n-2)<y^{-1}(n-1)$.
Then $y'\coloneqq s_1(1\ 2\ \cdots\ n)y(1\ 2\ \cdots\ n)^{-1}$ satisfies $(\ast)$ and $\ell(y')=\ell(y)+1$.
So by the induction hypothesis, there exist at least $j-1$ elements in $W_0$ which are $n$-cycles of length $3n-2j-5$ satisfying $(\ast)$.
Note that for any $r\in W_{\{s_2,\ldots,s_{n-3}\}}$, we have $r'y'r'^{-1}=s_1(1\ 2\ \cdots\ n)ryr^{-1}(1\ 2\ \cdots\ n)^{-1}$, where $r'=(1\ 2\ \cdots\ n)r(1\ 2\ \cdots\ n)^{-1}\in W_{\{s_2,\ldots,s_{n-2}\}}$.
So again by the induction hypothesis, it is easy to verify that there exists $r\in W_{\{s_2,\ldots,s_{n-3}\}}$ such that $r'y'r'^{-1}$ is a Coxeter element and $r'(\Phi_+\setminus \Phi_{\vp^\mu y'})\subset \Phi_+$.
Set $c=s_{n-2}s_{n-1}s_{n-3}\cdots s_2s_1$.
It is easy to check that if $n$ is odd (resp.\ even), then 
\begin{align*}
c,\  cs_{n-2}s_{n-3},\ \ldots\ ,\  cs_{n-2}s_{n-3}\cdots s_2,\  &cs_{n-2}s_{n-3}\cdots s_2 s_3s_4,\ \ldots\ ,\\
&cs_{n-2}s_{n-3}\cdots s_2 s_3s_4\cdots s_{n-2}s_{n-1}\\
(\text{resp.\ }c,\  cs_{n-2}s_{n-3},\ \ldots\ ,\  cs_{n-2}s_{n-3}\cdots s_3,\  &cs_{n-2}s_{n-3}\cdots s_3 s_2s_3,\ \ldots\ ,\\
&cs_{n-2}s_{n-3}\cdots s_3 s_2s_3\cdots s_{n-2}s_{n-1})
\end{align*}
are $n$-cycles satisfying $(\ast)$.
If $y'$ is one of the elements above, then $\Phi_{\{\chi_{2,3},\ldots,\chi_{n-2,n-1}\}}\cap\Phi_+\subset\Phi_{\vp^\mu y'}$ and there exists $r'\in W_{\{s_2,\ldots,s_{n-2}\}}$ such that $r'y'r'^{-1}$ is a Coxeter element.
Thus the claim is also true for $n$.
By induction, our claim is true for any $n\geq 4$.

Clearly $\nu_w=\nu_b$ for any $w\in \SAdm(\mu)_{\cyc}^\circ$.
Since $b=\tau^{n-1}$ is superbasic, the unique minimal length element in the $\sigma$-cojugacy class of $w$ is $\tau^{n-1}$ (cf.\ \cite[Proposition 3.5]{HN2}).
By Theorem \ref{minimal}, there exist a reduction tree $\cT$ for $w$ and a reduction path in $\cT$ such that $\en(\unp)=\tau^{n-1}$.
Thus by Lemma \ref{injective}, Proposition \ref{rational point}, Lemma \ref{semi-module 1n-2 lemm} and the claim we have shown above, there exist exactly $j$ elements of length $2j$ in $\SAdm(\mu)_{\cyc}^\circ$.
Moreover, it follows that $\pi(X_w(b)^0)$ is irreducible of dimension $\frac{\ell(w)}{2}$ for any $w\in \SAdm(\mu)_{\cyc}^\circ$ and that $X_w(b)=\emptyset$ for any $w\in \SAdm(\mu)\setminus \SAdm(\mu)_{\cyc}$.

It remains to show that the semi-module stratification of $X_{\mu}(b)$ is a refinement of the Ekedahl-Oort stratification.
We prove that for any $w\in \SAdm(\mu)_{\cyc}^\circ$, there exists an extended semi-module $(A^\ld,\vph)$ for $\mu$ such that $\pi(X_w(b)^0)=S_{A^\ld,\vph}(=X_\mu^\ld(b)$ by Lemma \ref{semi-module str} and Lemma \ref{semi-module 1n-2 lemm}).
We argue by induction on $\ell(w)$.
If $\ell(w)=2$, i.e., $w=\vp^\mu cs_{n-2}s_{n-3}\cdots s_2 s_3s_4\cdots s_{n-2}s_{n-1}=s_0s_{n-1}\tau^{n-1}$, then $w\rightarrow_{\sigma} s_0ws_0=\tau^{n-1}$.
It easily follows from Theorem \ref{empty} that $X_{\tau^{n-1}s_0}(b)=\emptyset$.
So by Proposition \ref{DL method prop}, we have $X_w(b)^0=Is_0I/I$ and hence $\pi(X_w(b)^0)=X_\mu^{\chi_{1,n}^\vee}(b)$.

Suppose that $\ell(w)\geq 4$ and the claim is true for any $w'\in \SAdm(\mu)_{\cyc}^\circ$ with $\ell(w')<\ell(w)$.
Since $\pi(X_w(b)^0)$ is irreducible of dimension $\frac{\ell(w)}{2}$, there exists a unique extended semi-module $(A^\ld,\vph)$ for $\mu$ such that $\dim (\pi(X_w(b)^0)\cap S_{A^\ld,\vph})=\frac{\ell(w)}{2}$.
Also, $\pi(X_w(b)^0)\cap S_{A^\ld,\vph}$ is open in both $\pi(X_w(b)^0)$ and $S_{A^\ld,\vph}$.
So the closure of $\pi(X_w(b)^0)\cap S_{A^\ld,\vph}$ in $X_\mu(b)$ is equal to both the closure of $\pi(X_w(b)^0)$ and $S_{A^\ld,\vph}$ in $X_\mu(b)$.
By \cite[Proposition 2.6]{He11} (see also \cite[\S 3.3]{GH}), the closure of $\pi(X_w(b)^0)$ is contained in $$\bigsqcup_{w'\in \SAdm(\mu)_{\cyc}^\circ, w'\le_S w}\pi(X_{w'}(b)).$$
Here we write $w'\le_S w$ if there exists $x\in W_0$ such that $xw'x^{-1}\le w$.
By the above description of the finite part of each element in $\SAdm(\mu)_{\cyc}^\circ$, it is easily checked that if $w'\in \SAdm(\mu)_{\cyc}^\circ$ and $\ell(w)=\ell(w')$, then there is no $x\in W_0$ such that $xwx^{-1}=w'$.
So if $w'\in \SAdm(\mu)_{\cyc}^\circ, w'\le_S w$ and $\ell(w')=\ell(w)$, then $w=w'$.
Thus by the induction hypothesis, we have $S_{A^\ld,\vph}\subseteq \pi(X_w(b)^0)$.
By \cite[Proposition 2.11 (5) \& Proposition 3.4]{CV}, the closure of $S_{A^\ld,\vph}$ is contained in a union of semi-module strata $T_\ld$ such that $\dim (T_\ld\setminus S_{A^\ld,\vph})<\dim S_{A^\ld,\vph}$.
Thus by the induction hypothesis and Lemma \ref{semi-module 1n-2 lemm}, we have $\pi(X_w(b)^0)\subseteq S_{A^\ld,\vph}$.
Therefore it follows that $\pi(X_w(b)^0)=S_{A^\ld,\vph}$, which completes the proof.
\end{proof}

\subsection{The Ekedahl-Oort Stratification for $\omega_1+\omega_{n-3}$}
\label{EO 1n-3}
Throughout this subsection, we set $\mu=\omega_1+\omega_{n-3}$.
Also we assume that $n\geq 7$.
Note that the unique dominant cocharacter $\mu'$ with $\mu'\prec \mu$ is $\omega_{n-2}$.

\begin{prop}
\label{EO 1n-3 not}
There exist at least $2(n-4)$ elements of length $3n-11$ in $\SAdm(\mu)_{\cyc}^\circ\coloneqq\SAdm(\mu)_{\cyc}\setminus \SAdm(\omega_{n-2})_{\cyc}$.
There also exists an element $w$ of length $3n-14$ in $\SAdm(\mu)$ such that its finite part is not $n$-cycle and $X_w(b)\neq \emptyset$.
Moreover, the semi-module stratification of $X_\mu(b)$ is not a refinement of the Ekedahl-Oort stratification.
\end{prop}
\begin{proof}
For any $1\le j\le n-4$, set $c_j=s_{n-3}s_{n-2}s_{n-1}s_{n-4}\cdots s_{j+2}s_{j+1}s_1\cdots s_{j-1}s_j$.
For $j=n-3$, set $c_{n-3}=s_1s_2\cdots s_{n-1}$.
Then we have $\vp^\mu c_j\in \SAdm(\mu)_{\cyc}^\circ$ and $\ell(\vp^\mu c_j)=3n-9$ for any $1\le j\le n-3$.
If $1\le j\le n-5$, then $c_js_{n-3}s_{n-2}$ and $c_js_{n-3}s_{n-4}$ are $n$-cycles of length $3n-11$ satisfying $\vp^\mu c_js_{n-3}s_{n-2},\vp^\mu c_js_{n-3}s_{n-4}\in \SAdm(\mu)_{\cyc}^\circ$.
Further $c_{n-4}s_{n-3}s_{n-2}$ and $c_{n-3}s_{n-4}s_{n-3}$ are also $n$-cycles of length $3n-11$ satisfying $\vp^\mu c_{n-4}s_{n-3}s_{n-2},\vp^\mu c_{n-3}s_{n-4}s_{n-3}\in \SAdm(\mu)_{\cyc}^\circ$.
Thus we have found $2(n-4)$ distinct elements of length $3n-11$ in $\SAdm(\mu)_{\cyc}^\circ$.

Set $y=c_{n-5}s_{n-3}s_{n-2}s_{n-4}s_{n-6}s_{n-5}=(1\ 2\ \cdots\ n-6\ n-2\ n\ n-3)(n-4\ n-5\ n-1)$.
Then $\vp^\mu y\in \SAdm(\mu)$ and $\chi_{1,n-1},\chi_{n-5,n}\in \Phi_+\setminus \Phi_{\vp^\mu y}$.
By Theorem \ref{empty}, $X_{\vp^\mu y}(b)\neq \emptyset$.
This shows the second assertion.
We can easily check the last assertion using Lemma \ref{semi-module 1n-3 lemm}, similarly as the proof of Proposition \ref{EO i not}.
\end{proof}

\subsection{The Ekedahl-Oort Stratification for $\omega_1+\omega_2,\omega_4+\omega_{n-1}$}
\label{EO 12}
Note that the unique dominant cocharacter $\mu'$ with $\mu'\prec \omega_1+\omega_2$ is $\omega_3$.
By an explicit calculation, it is easy to verify the following statements (cf.\ Proposition \ref{EO i refinement}).
\begin{prop}
\label{EO 12 refinement}
Assume that $n=5$.
Set $\mu=\omega_1+\omega_2$.
For any $1\le j\le 3(=\dim X_\mu(b))$, set $\SAdm(\mu)_{\cyc}^\circ\coloneqq\SAdm(\mu)_{\cyc}\setminus \SAdm(\omega_3)_{\cyc}$.
Then we have
\begin{align*}
\SAdm(\mu)_{\cyc}^\circ=\{s_0s_4s_3s_2s_1s_0\tau^3,s_0s_1s_4s_3s_0s_4\tau^3,s_0s_4s_3s_2\tau^3,s_0s_1s_4s_3\tau^3\}.
\end{align*}
Let $\vp^\mu y\in \SAdm(\mu)_{\cyc}^\circ$.
Then there exists $v\in \LP(\vp^\mu y)$ such that $v^{-1}yv$ is a Coxeter element.
Moreover, $X_w(b)=\emptyset$ for any $w\in \SAdm(\mu)\setminus \SAdm(\mu)_{\cyc}$, and the semi-module stratification of $X_{\mu}(b)$ is a refinement of the Ekedahl-Oort stratification.
\end{prop}

\begin{lemm}
\label{EO 12 not}
Assume that $n=7$ or $8$.
Let $\mu$ be $\omega_1+\omega_2$ (resp.\ $\omega_4+\omega_{n-1}$).
Set $c=s_1s_2\cdots s_{n-1}$.
Then $\vp^\mu cs_1s_2s_3\in \SAdm(\mu)$ and $X_{\vp^\mu cs_1s_2s_3}(b)\neq \emptyset$ (resp.\ $\vp^\mu c^{-1}s_5s_4s_3\in \SAdm(\mu)$ and $X_{\vp^\mu c^{-1}s_5s_4s_3}(b)\neq \emptyset$).
Further $cs_1s_2s_3$ (resp.\ $c^{-1}s_5s_4s_3$) is not $n$-cycle.
\end{lemm}

\subsection{The Ekedahl-Oort Stratification for $\omega_2+\omega_{n-3}$}
We set $\mu=\omega_2+\omega_{n-3}$.
Also we assume that $n\geq 5$.
\begin{lemm}
\label{EO 2n-3 not}
If $n$ is odd (resp.\ even), set $y=s_2s_3\cdots s_{n-3}s_1s_2\cdots s_{n-3}$ (resp.\ $y=s_2s_3\cdots s_{n-3}s_1s_2\cdots s_{n-2}$).
Then $\vp^\mu y\in \SAdm(\mu),X_{\vp^\mu y}(b)\neq \emptyset$ and $y$ is not $n$-cycle.
\end{lemm}
\begin{proof}
If $n$ is odd (resp.\ even), then $y=(1\ 3\ \cdots\ n-2)(2\ 4\ \ \cdots\ n-1\ n)$ (resp.\ $(1\ 3\ \cdots\ n-1)(2\ 4\ \cdots\ n)$) and $\vp^\mu y\in \SAdm(\mu)$.
Note that $\chi_{1,n},\chi_{2,n-1}\in \Phi_+\setminus \Phi_{\vp^\mu y}$.
So by Lemma \ref{empty}, $X_{\vp^\mu y}(b)\neq \emptyset$.
The proof is finished.
\end{proof}

\section{Comparison of Two Stratifications}
Keep the notations and assumptions in \S3.
\label{comparison of two stratifications}
\subsection{Known Cases}
The following results are known in (the proof of) \cite[Corollary 5.5 \& Theorem 5.9]{Shimada2}.
\begin{prop}
\label{known cases}
Let $\cong$ denote a universal homeomorphism.
\begin{enumerate}[(i)]
\item 
Assume that $n\geq 3$.
Set $\mu=2\omega_1, w=\vp^\mu s_1s_2\cdots s_{n-1}$ and
\begin{align*}
\ld=\begin{cases}
\chi_{2,n-1}^\vee+\chi_{4,n-3}^\vee+\cdots+\chi_{\frac{n-1}{2},\frac{n+3}{2}}^\vee & (\frac{n-1}{2}\ \text{even})\\
\chi_{1,n}^\vee+\chi_{3,n-2}^\vee+\cdots+\chi_{\frac{n-1}{2},\frac{n+3}{2}}^\vee & (\frac{n-1}{2}\ \text{odd}).
\end{cases}
\end{align*}
Then we have
$X_\mu(b)^0=X_\mu^\ld(b)=\pi(X_w(b)^0)\cong \A^{\frac{n-1}{2}}$.
\item 
Assume that $n\geq 3$.
Set $\mu=2\omega_1+\omega_{n-1},w_j=\vp^\mu s_{n-1}s_{n-2}\cdots s_{n-j+1}s_1s_2\cdots s_{n-j}$ and
\begin{align*}
\ld_j=\begin{cases}
\chi_{1,2j}^\vee+\chi_{2,2j-1}^\vee+\cdots+\chi_{j,j+1}^\vee & (j\le \frac{n}{2})\\
\chi_{2j+1-n,n}^\vee+\chi_{2j+2-n,n-1}^\vee+\cdots+\chi_{j,j+1}^\vee & (j\geq \frac{n}{2}).
\end{cases}
\end{align*}
for $j=1,2,\ldots,n-1$.
Then we have
$X_\mu(b)^0=\bigsqcup_{1\le j\le n-1} X_\mu^{\ld_j}(b)$ and $X_\mu^{\ld_j}(b)=\pi(X_{w_j}(b)^0)\cong \A^{n-1}$ for each $j$.
\item 
Assume that $n=5$.
Set $\mu=3\omega_1, w=\vp^\mu s_1s_2s_3s_4$ and $\ld=\chi_{1,2}^\vee+\chi_{3,4}^\vee$.
Then we have
$X_\mu(b)^0=X_\mu^\ld(b)=\pi(X_w(b)^0)\cong \A^4$.
\item 
Assume that $n=4$.
Set $\mu=3\omega_1, w=\vp^\mu s_1s_2s_3$ and $\ld=\chi_{3,2}^\vee$.
Then we have
$X_\mu(b)^0=X_\mu^\ld(b)=\pi(X_w(b)^0)\cong \A^3$.
\item 
Assume that $n=3$.
Set $\mu=4\omega_1, w=\vp^\mu s_1s_2$ and $\ld=\chi_{3,1}^\vee$.
Then we have
$X_\mu(b)^0=X_\mu^\ld(b)=\pi(X_w(b)^0)\cong\A^3$.
\item 
Assume that $n=3$.
Set $\mu=3\omega_1+\omega_2, w_1=\vp^\mu s_1s_2,w_2=\vp^\mu s_2s_1,\ld_1=\chi_{2,3}^\vee$ and $\ld_2=\chi_{3,2}^\vee$.
Then we have
$X_\mu(b)^0=X_\mu^{\ld_1}(b)\sqcup X_\mu^{\ld_2}(b)$ and $X_\mu^{\ld_j}(b)=\pi(X_{w_j}(b)^0)\cong \A^3$ for each $j$.
\item 
Assume that $n=2$.
Set $\mu=m\omega_1$ with $m\geq 1$, $w=\vp^\mu s_1$ and 
$$\ld=\begin{cases}
\frac{m-1}{2}\chi_{1,2}^\vee & (\frac{m-1}{2}\ \text{odd})\\
\frac{m-1}{2}\chi_{2,1}^\vee & (\frac{m-1}{2}\ \text{even}).
\end{cases}$$
Then we have
$X_\mu(b)^0=X_\mu^\ld(b)=\pi(X_w(b)^0)\cong\A^{\frac{m-1}{2}}$.
\end{enumerate}
\end{prop}

\subsection{Proof of the Main Theorem}
\begin{theo}
\label{main theo}
Let $\mu\in \Y_+$.
The following assertions are equivalent.
\begin{enumerate}[(i)]
\item The semi-module stratification of $X_{\pc\mu}(b)$ gives a refinement of the Ekedahl-Oort stratification.
\item For any $w\in {^S{\Adm}}(\mu)$ with $X_w(b)\neq \emptyset$, there exists $v\in \LP(w)$ such that $v^{-1}p(w)v$ is a Coxeter element.
\item The cocharacter $\mu$ has one of the following forms:
\begin{align*}
&\omega_1,\quad \omega_{n-1},\ &(n\geq 1),\\
&\omega_2,\quad 2\omega_1,\quad \omega_{n-2},\quad 2\omega_{n-1},\ &(\text{odd}\ n\geq 3),\\
&\omega_2+\omega_{n-1},\quad 2\omega_1+\omega_{n-1}\quad \omega_1+\omega_{n-2},\quad\omega_1+2\omega_{n-1},\ &(n\geq 3),\\
&\omega_3,\quad\omega_{n-3},\ &(n=7,8),\\
&3\omega_1,\quad 3\omega_{n-1},\ &(n=4,5),\\
&\omega_1+\omega_2,\quad\omega_3+\omega_4,\  &(n=5),\\
&4\omega_1,\quad \omega_1+3\omega_2,\quad 4\omega_2,\quad 3\omega_1+\omega_2, &(n=3),\\
&m\omega_1\ \text{with $m$ odd,} &(n=2).
\end{align*}
\end{enumerate}
If one of the above conditions holds, then for any $w\in \SAdm(\mu)_{\cyc}$, there exist $\mu'\in \Y_+$ with $\mu'\pc \mu$ and a cyclic extended semi-module $(A^\ld,\vph)$ for $\mu'$ such that $\pi(X_w(b)^0)=X_{\pc\mu}^\ld(b)=S_{A^\ld,\vph}$.
Moreover  $\pi(X_w(b)^0)$ is universally homeomorphic to an iterated fibration over a point whose fibers are (the perfection of) $\A^1$.
\end{theo}
\begin{proof}
For any $w=\vp^\mu y\in \SW$ with $\mu$ dominant, set $w^*=\vp^{(\mu(1),\ldots,\mu(1))}\varsigma(w)$ (cf.\ \S\ref{LP} and \S\ref{J-str semi-module}).
Then $w^*\in \SW$ and $p(w^*)=w_{\max}yw_{\max}^{-1}$ (cf.\ \S\ref{LP} and \S\ref{J-str semi-module}).
Note that the arguments and results in \S 5 and \S 6 for $(\mu,w,b)$ also hold for $(\mu^*,w^*,b^*)$.
Thus in this proof, it suffices to treat the case for either $\mu$ or $\mu^*$.

First assume that $n\geq 6$.
Let $1\le m_0<n$ be the residue of $m$ modulo $n$.
If $4\le m_0\le n-4$, then $\omega_{m_0}+\mn\omega_n\pc \mu$.
So by Lemma \ref{n=9} and Proposition \ref{EO i not}, $\mu$ satisfies neither (i) nor (ii).
If $n\geq 10$ and $m_0=3$, then by Lemma \ref{n=10}, $\mu$ satisfies neither (i) nor (ii).
If $n=7,8$ and $m_0=3$, then by Proposition \ref{EO i refinement}, $\mu=\omega_3$ satisfies (i) and (ii).
If moreover, $\mu\neq \omega_3$, then $\omega_1+\omega_2+\mn\omega_n\pc \mu$ or $\omega_4+\omega_{n-1}+(\mn-1)\omega_n\pc \mu$.
So by Lemma \ref{semi-module 12 n=7,8} and Lemma \ref{EO 12 not}, $\mu$ satisfies neither (i) nor (ii).
If $m_0=n-2$, then $\omega_1+\omega_{n-3}+\mn\omega_n\pc \mu$ unless $\mu=\omega_{n-2}$ or $2\omega_{n-1}$.
If $m_0=n-1$, then $\omega_2+\omega_{n-3}+\mn\omega_n\pc \mu$ unless $\mu=\omega_{n-1},\omega_1+\omega_{n-2}$ or $\omega_1+2\omega_{n-1}$.
Thus the equivalence of (i), (ii) and (iii) for $m_0=n-2,n-1$ follows from Theorem \ref{cyclic classification}, Proposition \ref{EO i refinement}, Proposition \ref{EO 1n-3 not}, Proposition \ref{EO 2n-3 not} and Proposition \ref{known cases}.

Assume that $n=5$.
If $m_0=3$, then $\omega_1+\omega_3+\omega_4+\mn\omega_n\pc \mu$ unless $\mu=\omega_{3},2\omega_4,\omega_1+\omega_2$ or $3\omega_1$.
If $m_0=4$, then $2\omega_2+\mn\omega_n\pc \mu$ unless $\mu=\omega_{4},\omega_1+\omega_{3}$ or $\omega_1+2\omega_{4}$.
Set $y_5=(1\ 5\ 3)(2\ 4)$.
Then it is easy to check that $\vp^{\omega_1+\omega_3+\omega_4} y_5\in \SAdm(\omega_1+\omega_3+\omega_4)$ and $X_{\vp^{\omega_1+\omega_3+\omega_4} y_5}(\tau^8)\neq \emptyset$.
Assume that $n=4$.
If $m_0=3$, then $2\omega_2+\omega_3+\mn\omega_n\pc \mu$ unless $\mu=\omega_{3},\omega_1+\omega_2,\omega_1+2\omega_3$ or $3\omega_1$.
Set $y_4=(1\ 3)(2\ 4)$.
Then it is easy to check that $\vp^{2\omega_2+\omega_3} y_4\in \SAdm(2\omega_2+\omega_3)$ and $X_{\vp^{2\omega_2+\omega_3} y_4}(\tau^7)\neq \emptyset$.
Assume that $n=3$.
If $m_0=2$, then $2\omega_1+3\omega_2+\mn\omega_n\pc \mu$ unless $\mu=\omega_{2},2\omega_1,\omega_1+2\omega_2,3\omega_1+\omega_2$ or $4\omega_2$.
Set $y_3=(1\ 3)$.
Then it is easy to check that $\vp^{2\omega_1+3\omega_2} y_3\in \SAdm(2\omega_1+3\omega_2)$ and $X_{\vp^{2\omega_1+3\omega_2} y_3}(\tau^8)\neq \emptyset$.
Thus the equivalence of (i), (ii) and (iii) for $n=2,3,4,5$ also follows from Theorem \ref{cyclic classification}, Proposition \ref{EO i refinement}, Proposition \ref{EO 2n-3 not} and Proposition \ref{known cases}.
The case for $n=1$ is trivially true.

Assume that $\mu$ satisfies one of the conditions in the theorem, which is equivalent to each other as we have just proved.
Then the last assertion follows from the proof of the results in \S5 and \S6 except the ``moreover'' part.
Let $w\in \SAdm(\mu)_{\cyc}$.
By Lemma \ref{injective}, the map $X_w(b)^0\rightarrow \pi(X_w(b)^0)$ induced by $\pi$ is universally bijective.
Using \cite[Proposition 3.1.1]{GH} and Proposition \ref{DL method}, we can easily check that $\pi^{-1}(\pi(X_w(b)^0))\cap (\bigcup_{w'\le w}X_w(b))=X_w(b)^0$ (cf.\ the proof of \cite[Lemma 5.8]{Shimada2}).
Since $\pi$ is proper, the map $X_w(b)^0\rightarrow \pi(X_w(b)^0)$ is also proper.
This implies that this map is a universally homeomorphism.
Again by the results in \S5 and \S6, $X_w(b)^0$ and hence $\pi(X_w(b)^0)$ are universally homeomorphic to an iterated fibration over a point whose fibers are (the perfection of) $\A^1$.
This finishes the proof.
\end{proof}

\begin{rema}
\label{main rema}
Except the cases where $\mu$ or $\mu^*$ is $\omega_1+\omega_{n-2}\ (n\geq 3)$ or $\omega_1+\omega_2\ (n=5)$, it follows from \cite[Theorem 5.3]{Viehmann2} and Theorem \ref{known cases} that if one of the conditions in Theorem \ref{main theo} holds, then each $X_{\mu}^\ld(b)(\neq \emptyset)$ is universally homeomorphic to an affine space.
We do not know this is true in general.
\end{rema}

\bibliographystyle{myamsplain}
\bibliography{reference}
\end{document}